%% file: main-IJNBME-r2.tex
\documentclass{article}
 \usepackage{a4wide}

\usepackage{amsfonts}
\usepackage{graphicx}  
\usepackage{hyperref}
\usepackage{pdfsync}
\usepackage{ulem}
\graphicspath{{./figures/}}
\usepackage{colortbl}
\usepackage{pgfplots}
\usepackage{pgfplotstable}
\usepackage{tikz}
\usepackage{tkz-euclide}
\usetkzobj{all}
\usetikzlibrary{decorations.pathreplacing,shapes.misc,calc,intersections,arrows,external}
\usepgfplotslibrary{external,statistics}
\usepackage{xifthen}
\usepackage{stmaryrd}
\usepackage{definitions}

\newcommand{\tens}[1]{\uline{\boldsymbol #1}}

\newcommand{\btau}{\boldsymbol{\tau}}

\begin{document}

\title{Multiscale modeling of vascularized tissues via non-matching immersed methods}

\author{Luca Heltai\footnote{International School for Advanced
    Studies, Trieste, Italy} \and Alfonso Caiazzo\footnote{Weierstrass
    Institute for Applied Analysis and Stochastics, Berlin,
    Germany}\phantom{.}\footnote{Corresponding author: \texttt{caiazzo@wias-berlin.de}}}

\maketitle





\section*{Abstract}

We consider a multiscale approach based on immersed methods for the
  efficient computational modeling of tissues composed of an elastic
  matrix (in two or three-dimensions) and a thin vascular structure
  (treated as a co-dimension two manifold) at a given pressure.  We
  derive different variational formulations of the coupled problem, in
  which the effect of the vasculature can be surrogated in the
  elasticity equations via singular or hyper-singular forcing
  terms. These terms only depend on information defined on
  co-dimension two manifolds (such as vessel center line, cross
  sectional area, and mean pressure over cross section), thus
  drastically reducing the complexity of the computational model.  We
  perform several numerical tests, ranging from simple cases with
  known exact solutions to the modeling of materials with random
  distributions of vessels. In the latter case, we use our immersed
  method to perform an in silico characterization of the mechanical
  properties of the effective biphasic material tissue via statistical
  simulations.


\section{Introduction}
\label{sec:introduction}
This paper is motivated by mathematical and computational modeling in
the context of tissue imaging, such as Magnetic Resonance Elastography
(MRE), a quantitative imaging technique sensitive to the mechanical
properties of living tissues.
In MRE, the tissue undergoes external harmonic excitations, such as
shear or compression waves, and Magnetic Resonance Imaging is used to
recover the mechanical response of the tissue in terms of the internal
displacement.  These data, combined with a mathematical model of the
underlying tissue dynamics, are then employed to characterize the
tissue -- in vivo and non invasively

The \textit{inversion process}, 
i.e., the recovery of mechanical property from displacement data, is mainly based
on simplified models, assuming that tissues behaves as homogeneous and isotropic linear elastic or viscoelastic
materials. Such models have been already used to demonstrate the potential of MRE in characterizing
pathological tissues (e.g., cancer of fibrosis) \cite{muthupillai-96,sack-2008,sack-book,wuerfel-10}.

However, in several clinical applications, the complex multiphysics and multiscale nature of living tissues cannot be
neglected. The characterization of vascularized tissues is one of these examples.

Recent experimental studies have been dedicated to understanding
the potential of
MRE to characterize intrinsic properties of biphasic tissues (e.g.,
brain and liver), aiming at the non invasive diagnosis of pressure
related diseases (see, e.g.,
\cite{hirsch-etal-2013-compmre,hirsch-etal-2014-liver}). 
Experiments comparing parameter estimation in vivo and ex vivo 
 confirmed that inversion methods deliver very different results
if the vascular component is inactive (see, e.g. \cite{chatelin-2011}) and/or if
the vasculature pressure varies \cite{hirsch-etal-2013-compmre}, thus requiring the 
development of more detailed models,
able to describe both phases (solid matrix and fluid
vasculature) are necessary.
On the other hand, in order to obtain results in a clinically relevant time, 
mathematical and computational shall be able to efficiently deal with 
the multiscale structure of the system.

From the computational point of view, fully resolved biphasic models,
i.e., accounting for the coupling between the tissue and fluid
vasculature at the microscale, are practically unfeasible. On the one
hand, the high geometrical complexity would lead to excessive
computational cost and, on the other hand, image resolution used in
MRE (of the order of millimeters) does not allow to reconstruct in
full detail the vasculature.
%

The goal of this work is to propose and test a novel mathematical
multiscale model for vascularized tissues (composed of an elastic
matrix and a thin -- pressurized -- fluid vasculature) with the purpose
of providing an efficient effective material model to be used for tissue characterization.
In the presented approach, 
the vasculature
(microscale) is not explicitly discretized, but it is immersed in the
elasticity problem describing the matrix dynamics at the macroscale.
To this aim, we use an approach based on the
Immersed Boundary Method (IBM), in order to account for complex (one-dimensional) structures
within two- and three-dimensional elastic materials.

The Immersed Boundary Method was introduced by Peskin
in~\cite{Peskin1972}, to study the blood flow around heart valves (see
also~\cite{Peskin2002}, or the review~\cite{Mittal2005b}), and evolved
into a large variety of methods and algorithms. The main idea behind
this technique is to address complex fluid-structure interaction
problems by formulating them on a fixed background fluid problem, with
the addition of singular source terms that take into account the
presence of the solid equations, removing the requirement that the
position of the interfaces between the fluid and solid domains should
be aligned with the computational mesh. 

In the original Immersed Boundary Method~\cite{Peskin1972} the
singular source terms are formally written in terms of the Dirac delta
distribution, and their discretization follows two possible routes: i)
the Dirac delta distribution is approximated through a smooth
function, or ii) the variational definition of the Dirac distribution
is used directly in the Finite Element formulation of the problem. For
finite difference schemes, the first solution is the only viable
option, even though the use of smooth kernels may excessively smear
the singularities, leading to large errors in the
approximation~\cite{Hosseini2014}. In the context of finite elements,
both solutions are possible. The methods derived from the
Immersed Finite Element Method (IFEM) still use approximations of the
Dirac delta distribution through the Reproducing Kernel Particle
Method (RKPM)~\cite{ZhangGerstenberger-2004-Immersed-finite-0}.

Variational formulations of the IBM were introduced
in~\cite{BoffiGastaldi-2003-a,
  BoffiGastaldi-2007-Numerical-stability-0,
  BoffiGastaldiHeltaiPeskin-2008-a, Heltai-2008-a}, and later
generalised in~\cite{Heltai2012b} and~\cite{RoyHeltaiCostanzo-2015-a},
where the need to approximate Dirac delta distributions is removed by
exploiting directly the weak formulation. Such formulations allow the
solution of PDEs with jumps in the gradients without enriching the
finite element space, and without introducing approximations of the
Dirac delta distribution.
In the context of 3D-1D multiscale models, an approach using
techniques similar to the IBM has been described in
\cite{dAngelo-08,dAngelo-12} for the case of diffusion equations. In
this case, a diffusion problem was solved on both the 3D and on the 1D
domains, considering, additionally, the 1D vasculature as a source of
nutrients for the 3D tissue.  This approach has been recently extended
to the case of a 3D porous media (Darcy) coupled to an immersed
vasculature, resolving the flow in the vascular network via a (0D)
lumped parameter model~\cite{cattaneo-14}.

In this paper, we consider the case of a 3D (or 2D) 
elastic matrix with an immersed 1D (resp. 0D) vasculature with a given fluid
pressure, i.e., under the assumption that the diameter of the fluid vessel
is much smaller than the size of the characteristic domain.

In the variational formulation, the effect of the fluid is then
included in the elasticity equations by means of a singular source
term on a lower dimensional manifold.  We begin by analyzing a
singular formulation in which the source term is concentrated on the
vessel boundary. Next, we discuss a \textit{hyper-singular}
alternative, in which the immersed source term is applied only at the
vessel centerline (a co-dimension two manifold), thus reducing
drastically the computational effort.

The multiscale model will be derived starting from a 2D-0D
axis-symmetric case and subsequently extended to the general 3D-1D
situation.  We perform different numerical tests, validating the model
in a simple setting in which an analytical solution is available, and
investigating the statistical effective behavior of a biphasic
material with random vessel distribution as a function of elastic and
geometrical parameters.  We focus on the effective tissue dynamics
assuming a steady known fluid pressure in the vasculature.

An extension of this model including a two-way coupling with an active
one-dimensional vasculature (e.g., using the approach described in
\cite{mueller-etal-16}) is currently under investigation and will be
subject of a future work.
%

The rest of the paper is organized as follows. In Section \ref{sec:2d}
we discuss the two-dimensional case, starting from a model problem
with known exact solution. 
The approach is extended in Section
\ref{sec:3d} to three dimensions.
In Section \ref{ssec:tissue-homo} 
we discuss a homogenized model for a pressurized tissue based on the singular formulation, and its implication
concerning the in silico characterization of mechanical properties.
The discretized model is described
in Section \ref{ssec:discrete}, while numerical results are
presented in Section \ref{sec:numerics}.Finally, Section
\ref{sec:conclusions} draws conclusions and future directions of our work.

\section{The two-dimensional model}\label{sec:2d}

\subsection{A simple problem setting}\label{ssec:2d-setting}
We consider the situation of a biphasic tissue composed of an
elastic matrix and thin blood vessels, under the assumption that the
vessel diameters are much smaller of the typical size of the
surrounding matrix. 
To fix the ideas, we start with the derivation for a two-dimensional
model problem, considering a single vessel.  Assuming that the vessels
are small compared to the elastic matrix, and that long term
interaction can be neglected, the arguments can be extended also to
general domains and multiple vessels.

\begin{figure}[!h]
\centering
\includegraphics{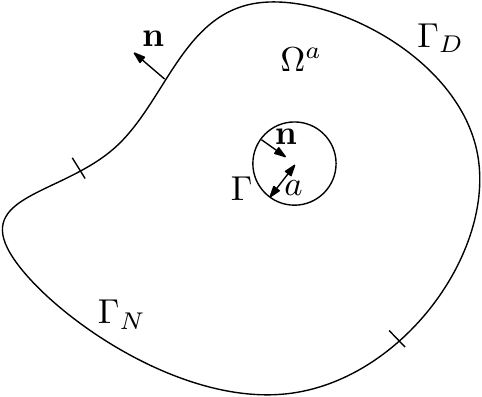}
\caption{An example of a domain $\Omega$ with a single vessel of radius $a$.}
\label{fig:2d-domain}
\end{figure}

Let $a >0$, and let us introduce the set
\[
\Ve{a} = \left\{ \bx \mid \| \bx \|  \leq a \right\}
\]
describing a circle of radius $a$ (which will be also referred to as
\textit{vessel}).  Next, let $\Omega \subset \mathbb R^2$ and let us
introduce the \textit{tissue} domain
$\Omega^a = \Omega \backslash \Ve{a}$. We assume that the boundary of
$\Omega$ is decomposed as
\[
\partial \Omega = \Gamma_D \cup \Gamma_N\,,
\]
we define $\Gamma := \partial \Ve{a}$ to be the vessel boundary,
and we denote with $\mathbf n$ the normal vector to $\Gamma$ pointing
outwards the tissue domain (see, e.g., the sketches in Figure
\ref{fig:2d-domain}).

This setting represents the case of a tissue that extends indefinitely
along the $z$-direction with an embedded cylindrical vessel with
cross-section equal to $\Ve{a}$.
We assume that the presence of a flow inside the vessel can be modeled
as a constant \emph{excess pressure} $\bar p$, that represents the
difference in pressure between the interior part of the vessel, and
the surrounding elastic matrix. The excess pressure represents the
force per unit area that the vessel exerts on the elastic matrix, and
we assume that this force is directed along the normal to the vessel.

We now consider the following problem:
\begin{problem}[2D on $\Omega^a$] \label{pb:two-dimensional} 
Given the excess pressure $\overline{p} > 0 $, find the displacement 
$\bue: \Omega^a \to \mathbb R^2$ solution to:
  \begin{equation}
    \label{elasto-omeps}
    \begin{aligned}
 & -\nabla \cdot \tens{\sigma}(\bue) = \boldsymbol 0,\;\text{ in }\;\Omega^a\\
& \bue = \boldsymbol 0,\; \text{ on }\; \Gamma_D \\
& \tens{\sigma}(\bue) \cdot \bn = \boldsymbol 0,\; \text{ on }\; \Gamma_N\\
& \tens{\sigma}(\bue) \cdot \bn = 
-\overline{p} \, \bn \; \text{ on }\; \Gamma
    \end{aligned}
  \end{equation}
The above system of equations describes the dynamics of 
a compressible, linear elastic material,   where 
\begin{equation}\label{eq:sigma-u}
\tens{\sigma}(\bu) := 2 \mu \tens{e}(\bu) + \lambda I \ldiv\bu,
\end{equation}
stands for the Cauchy  
stress tensor,
$\tens{e}(\bu) =\frac12\( \nabla \bu + \nabla \bu^T\)$ 
denotes the symmetric
part of the infinitesimal strain tensor, $\mu$ and $\lambda$ are the so
called Lam\'e constants, and $I$ is the identity matrix.
\end{problem}

\begin{remark}
  Notice that, using \eqref{eq:sigma-u}, the normal component of the
  solid stress $\tens{\sigma}(\bu) \cdot \bn$ can be also written as
  $\tens{\sigma}(\bu) \cdot \bn = (2\mu + \lambda) (\ldiv\bu) \bn$
\end{remark}

Let us now introduce the functional spaces
\begin{equation}
  \label{eq:definition-V}
    V^a := \{ \vv \in (H^1(\Omega^a))^2, \text{ such that }
    \vv|_{\Gamma_D} = \bs{0} \},
\end{equation}
Multiplying \eqref{elasto-omeps} with $\vv \in V^a$ and integrating by
parts yields a standard variational formulation of
Problem~\ref{pb:two-dimensional}:
\begin{problem}[2D on $\Omega^a$, variational] 
  \label{pb:two-dimensional-variational}
  Given the excess pressure $\overline{p} > 0 $, find the displacement
  $\bue \in V^a$ solution to:
  \begin{equation}
    \label{elasto-omeps-variational}
    \begin{aligned}
      & 2\mu(\tens{e}(\bue), \tens{e}(\vv))_{\Omega^a} + \lambda(\ldiv\bue, \ldiv
      \vv)_{\Omega^a} = \int_{\Gamma} \bar p \bn\cdot\vv \d \Gamma &&
      \forall \vv \in V^a,
    \end{aligned}
  \end{equation}
  where $(\cdot, \cdot)_{\Omega^a}$ denotes the inner product in
  $(L^2(\Omega^a))^2$.
\end{problem}

\subsection{Exact solution in the axis-symmetric
  case}\label{ssec:2d-exact}
\begin{figure}[!h]
\centering
\includegraphics{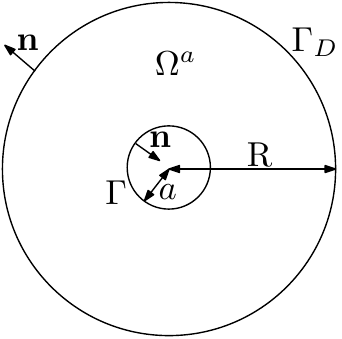}
\caption{An example of an axi-symmetric domain $\Omega$ with a single vessel of radius $a$.}
\label{fig:2d-domain-axi-symmetric}
\end{figure}

\begin{figure}[htp]
  \begin{center}
    \includegraphics[width=.45\textwidth]{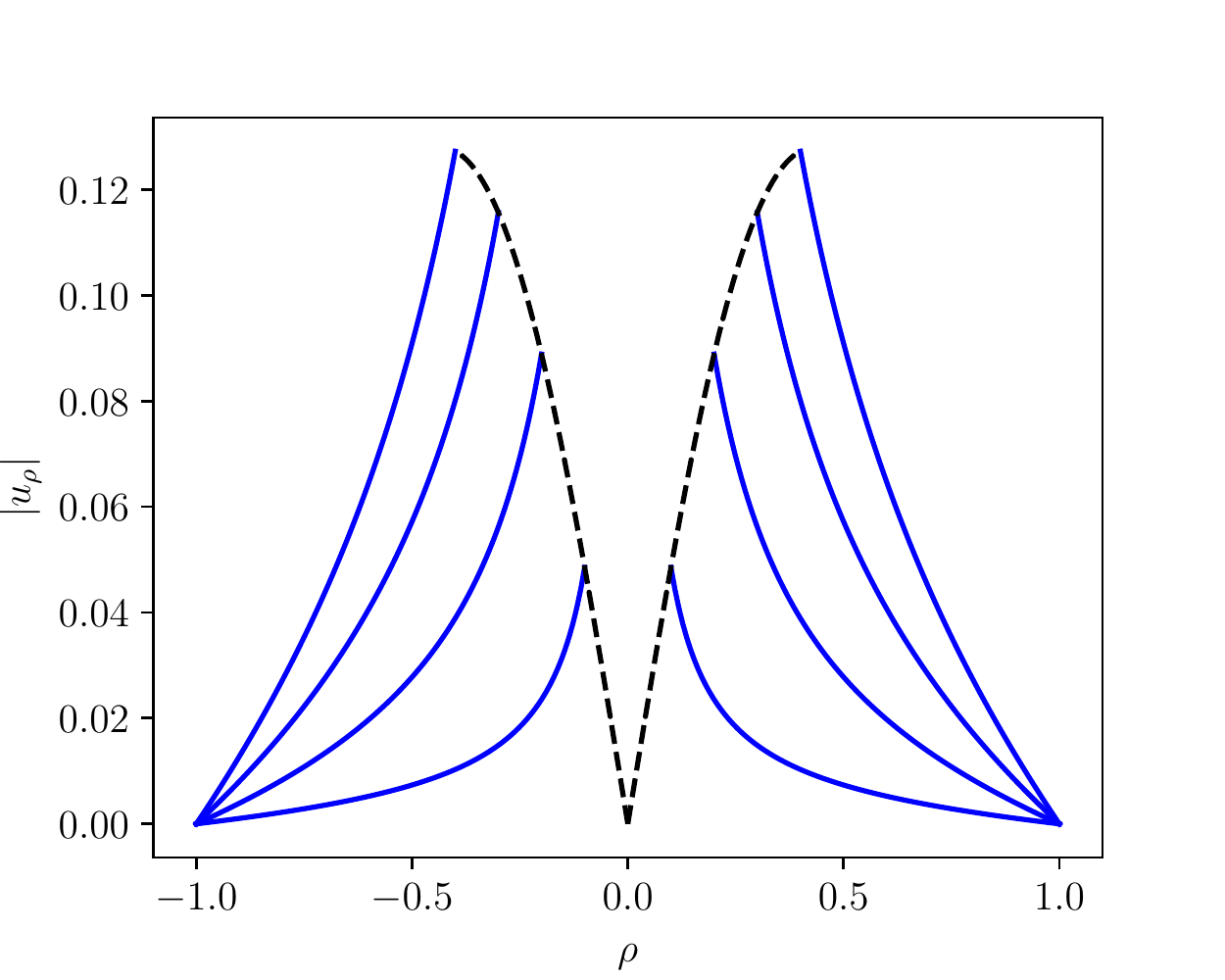}\hfill 
        \includegraphics[width=.45\textwidth]{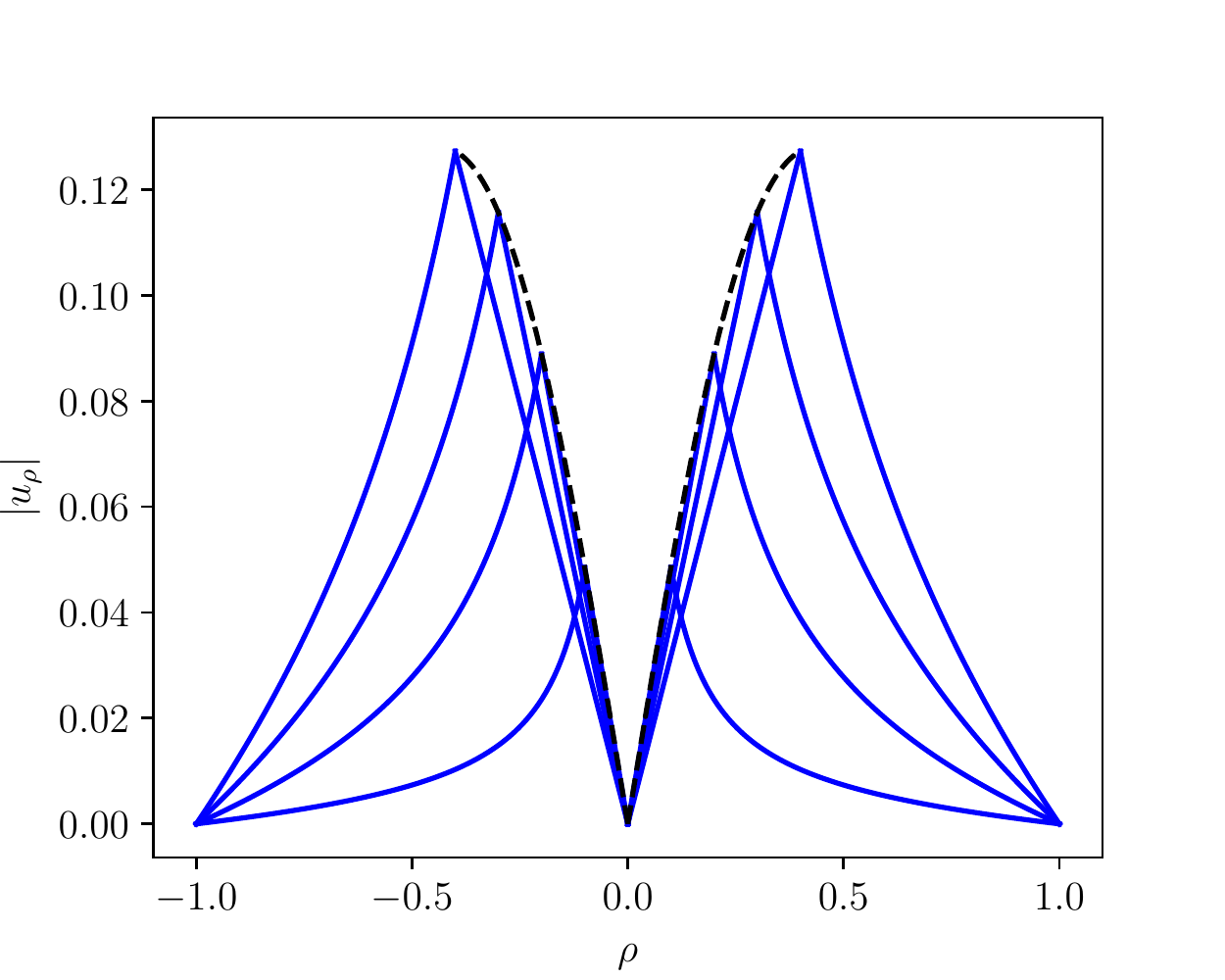}
  \end{center}
  \caption{Left: Central sections of the radial component of the exact
    solution $\bu_a$, for different values of the vessel radius
    $a \in [0.02,0.55]$ (with $\lambda=\mu=1$, $R=1$, $p=1$).
    Right: Central sections of the radial component of the exact
      solution $|u_\rho|$ to
      Problem~\ref{pb:two-dimensional-distributional}, for different
      values of the vessel size $a=.1, .2, .3, .4$ ($\lambda=\mu=1$, $R=1$,
      $p=1$).} \label{fig:exact-2d}
\end{figure}

In the special case where $\Omega$ is a circle of radius $R > a$ 
(see Figure~\ref{fig:2d-domain-axi-symmetric}), problem \ref{elasto-omeps}
can be solved analytically.  
In this case, due to the radial symmetry of the domain,
the angular component of the
solution vanishes, while the radial component depends only on the
distance from the vessel, i.e.,
$\bue (\bx) = u_\rho \frac{\bx}{|\bx|}$.  The elasticity problem
reduces therefore to an ODE for $u_\rho$, yielding the displacement
\begin{equation}\label{eq:2d-exact}
      \bue (\bx) = 
      \frac{\overline{p} a^{2} \left(R^{2} -
          |\bx|^2\right)}{2 \left(R^{2} \mu + \lambda a^{2} + 
          \mu a^{2}\right)} \frac{\bx}{|\bx|^2}.
    \end{equation}
Figure~\ref{fig:exact-2d}(left) shows the behaviour of the radial displacement given by 
\eqref{eq:2d-exact}, varying the size $a$ of the vessel. 

We are interested in the situation where the radius of the vessel is
small compared to the size of the tissue domain.  Under the assumption
$a \ll R$, we obtain the approximation:
 \begin{equation}\label{eq:2d-exact-bd}
   \bue |_{\Gamma} = \frac{\overline{p} a \left(R^{2} -
       a^{2}\right)}{2 R^{2} \mu + 2 \lambda a^{2} + 2
     \mu a^{2}} \frac{\bx}{|\bx|} = -\frac {\overline{p} a}{2\mu} \bn + O\left( \left(\frac{a}{R}\right)^2 \right)\,
    \end{equation}
where we also used the fact that, on $\Gamma$, it holds $\dfrac{\bx}{|\bx|} = - \bn$.

From \eqref{eq:2d-exact}-\eqref{eq:2d-exact-bd}, we conclude that the 
excess pressure  $\overline{p}$  induces a radial deformation of the surrounding elastic matrix 
(normal to the vessel boundary) which is of the order of $\frac {\overline{p} a}{2\mu}$ on the vessel boundary
and decays as $ \frac{1}{|\bx|}$.

\subsection{A singular problem on the whole domain}\label{ssec:2d-sing}
Next, we aim at formulating an extension of the axi-symmetric problem on the whole
domain $\Omega$, and at introducing a forcing term so that the solution of the extended problems
coincides with  the solution $\bue$, defined in \eqref{eq:2d-exact}, only outside of $\Ve{a}$.

In practice, we first extend the solution $\bue$ inside $\Ve{a}$ as a uniform deformation, i.e., linearly
in the distance from the origin (see Figure \ref{fig:exact-2d}, right):
\begin{equation}\label{eq:2d-extended}
\bu^{\Omega}(\bx) = \begin{cases}
 \dfrac{\overline{p}  \left(R^{2} -
          |\bx|^2\right)}{2 \left(R^{2} \mu + \lambda a^{2} + 
          \mu a^{2}\right)}\,  \dfrac{a^2 \bx}{|\bx|^2} & |\bx| \geq a \\[2em]
\dfrac{\overline{p} \left(R^{2} -
          a^{2}\right)}{2 R^{2} \mu + 2 \lambda a^{2} + 2
        \mu a^{2}} \,{\bx}   & |\bx| < a
\end{cases}
\end{equation}

The function $\bu^{\Omega}$ defined in \eqref{eq:2d-extended}
is continuous across the vessel boundary $\Gamma$.
However, the normal stress has a jump given by
\begin{equation}\label{eq:2d-sigma-jump}
  g_a  \bn  :=  \llbracket \tens{\sigma}(\bu^{\Omega}) \bn\rrbracket_{\Gamma} := \frac{R^{2} \overline{p}
        \left(\lambda + 2 \mu\right)}{R^{2} \mu + \lambda a^{2}
        + \mu a^{2}} \bn \,.
\end{equation}
Hence, in order to define an elasticity problem on $\Omega$,
whose solution is given by $\bu^{\Omega}$, we will consider
a fictitious elastic material defined on the whole domain, with the same properties as the original one
(defined in $\Omega^a$), but subjected to
a singular source term that imposes the jump $g_a  \bn$ in the normal stress.
Namely, we consider the following problem:
\begin{problem}[2D, singular] \label{pb:two-dimensional-distributional} 
Given an excess pressure $\overline{p}>0$,  find the distributional solution $\bu$ to:
  \begin{equation}
    \label{eq:two-d-distributional}
        \begin{aligned}
 & -\nabla \cdot \tens{\sigma}(\bu) = \bF^S_a,\;& &\text{ in }\;\Omega\\
& \bu = \boldsymbol 0,\; &&\text{ on }\; \Gamma_D \\
& \tens{\sigma}(\bu) \cdot \bn =  \boldsymbol 0,\; &&\text{ on }\; \Gamma_N
    \end{aligned}
    \end{equation}
 with
\begin{equation}\label{F2D-sing}
\bF^S_a(\bx) :=  \int_{\Gamma} \delta(\bx-\by) g_a \bn(\by)
  \d \Gamma_{\by},\; \qquad \forall  \bx \in \Omega
\end{equation}
where $\delta$ denotes the two-dimensional Dirac delta distribution
and $\mathbf y$ stands for a local coordinate on the interface $\Gamma$.
\end{problem}

In order to understand the definition of
Problem~\ref{pb:two-dimensional-distributional}, let us introduce the
Sobolev space 
\[
 V := \{ \vv \in (H^1(\Omega))^2, \text{ such that }
   \vv|_{\Gamma_D} = \bs{0} \}\,,
\]
 denoting with
$(\cdot,\cdot)$ the scalar product in $L^2(\Omega)$, and with $<\cdot, \cdot>$
the duality product between $H_0^1(\Omega)$ and its dual space
$H^{-1}(\Omega)$.  Multiplying \eqref{eq:two-d-distributional} with a
function $\vv \in V$ and integrating by parts over $\Omega^a$ we obtain:
 \begin{equation}\label{elasto-omeps-weak}
\begin{aligned}
 2 \mu \( \tens{e}(\bu),\tens{e}(\vv)\)_{\Omega^a} + \lambda \(\ldiv\bu,\ldiv\vv\)_{\Omega^a}   
 - \( \tens{\sigma}(\bu) \cdot \bn,\vv\)_{\Gamma} = 0\,.
\end{aligned}
\end{equation}
Proceeding similarly, but considering a fictitious elasticity problem
inside $\Ve{a}$ with the same characteristics of the surrounding
elastic matrix, we obtain 
\begin{equation}\label{elasto-Ve-weak}
      2\mu \( \tens{e}(\bu), \tens{e}( \vv ) \)_{\Ve{a}} + 
      \lambda \(\ldiv\bu,\ldiv\vv\)_{\Ve{a}} +
       \( \tens{\sigma}(\bu) \cdot \bn, \vv \)_{\Gamma} = 0
    \end{equation}
where the signs of the last terms in \eqref{elasto-omeps-weak} and \eqref{elasto-Ve-weak} 
depend on the chosen
 orientation of the normal vector $\bn$ (from the tissue towards the vessel). 
 Summing  \eqref{elasto-omeps-weak} and \eqref{elasto-Ve-weak},
 imposing continuity on the displacement and the given jump of the normal stress 
\eqref{eq:2d-sigma-jump}, we obtain the weak formulation:
\begin{problem}[2D, singular, variational]
  \label{pb:two-dimensional-singular-variational}
  Given an excess pressure $\bar p$, find the solution $\bu \in V$ such
  that 
  \begin{equation}\label{elasto-Omega-weak}
    (2\mu \tens{e}(\bu), \tens{e}(\vv))_{\Omega} + (\lambda \ldiv\bu, \ldiv
    \vv)_{\Omega} =  
    \int_{\Gamma} g_a \bn \cdot \vv \, \qquad \forall \vv \in V.
  \end{equation}
\end{problem}
Now, let us introduce the distributional definition of the two dimensional Dirac delta
distribution, i.e.,        
\begin{equation}
  \label{eq:Dirac-definition}
  \int_\Omega  \vv(\bx) \delta(\bx-\by) \d x = \vv(\by) \quad \forall \vv\in V
  \cap C^0(\Omega), \forall \by \in \Omega\,.
\end{equation}
Using \eqref{eq:Dirac-definition}, switching the order of integration,
and interpreting the integral on $\Gamma$ of functions in $V$ in
the sense of traces, it is possible to rewrite formally the term
$\int_\Gamma g_a \bn \cdot \vv$ as 
\begin{equation}
 \begin{split}
\int_{\Gamma} g_a(\by) \bn(\by) \cdot \phantom{\int_\Omega}  \vv(\by)
\phantom{\delta(\bx-\by) \d x} \d \Gamma_{\mathbf y} &
= \\
\int_{\Gamma} g_a(\by) \bn(\by) \cdot \int_\Omega  \vv(\bx) \delta(\bx-\by) \d x  \d \Gamma_{\mathbf y} &= \\
\int_{\Omega} ~\int_\Gamma g_a \bn(\mathbf y)
\delta(\bx-\mathbf y) \d \Gamma_{\mathbf y} ~\cdot~\vv(\bx)  \d \bx & =:
<\bF_a^S,\vv>,
\end{split}
\end{equation}
where $\bF_a^S$ is the singular forcing term introduced in \eqref{F2D-sing}.
%
%
%

For a detailed discussion on the behaviour of this distributional
forcing term, see~\cite{HeltaiRotundo-2016-a}. The term $\bF_a^S$ is a
distribution in $H^{-1}(\Omega)$, and was introduced originally
in~\cite{BoffiGastaldi-2003-b} and later generalized
in~\cite{Heltai-2008-a, BoffiGastaldiHeltaiPeskin-2008-a,
  HeltaiCostanzo-2012-a} as a variational formulation of the Immersed
Boundary Method~\cite{Peskin-2002-The-immersed-boundary-0}, to
approximate fluid structure interaction problems using non-matching
grids between the immersed structure and the surrounding fluid.

\subsection{The hypersingular problem}\label{ssec:hyper-2d}
The variational formulation introduced in \eqref{elasto-Omega-weak}
allows to reformulate the coupled problem as an elasticity problem on the whole domain $\Omega$
without explicitly taking into account the boundary condition on the
vessel boundary, so that, in the axi-symmetric case, the solution
coincides with the exact one outside the vessel $\Ve{a}$.
 From the practical point of view, this approach might be used to employ a spatial discretization (mesh) that
 does not explicitly resolve in full detail the vessel boundary, hence considerably reducing
the overall complexity, especially in the case of thin vessels.
 On the other hand, a discretization of the formulation described in
 Problem~\ref{pb:two-dimensional-distributional} still requires a
 characteristic mesh size that resolves the vessels boundary, in order 
 to compute accurately enough the integral on $\Gamma_a$. For small
 to very small vessels sizes, this constrain might still yield an excessive computational
 cost.

\begin{figure}
   \centering
  \includegraphics[width=.6\textwidth]{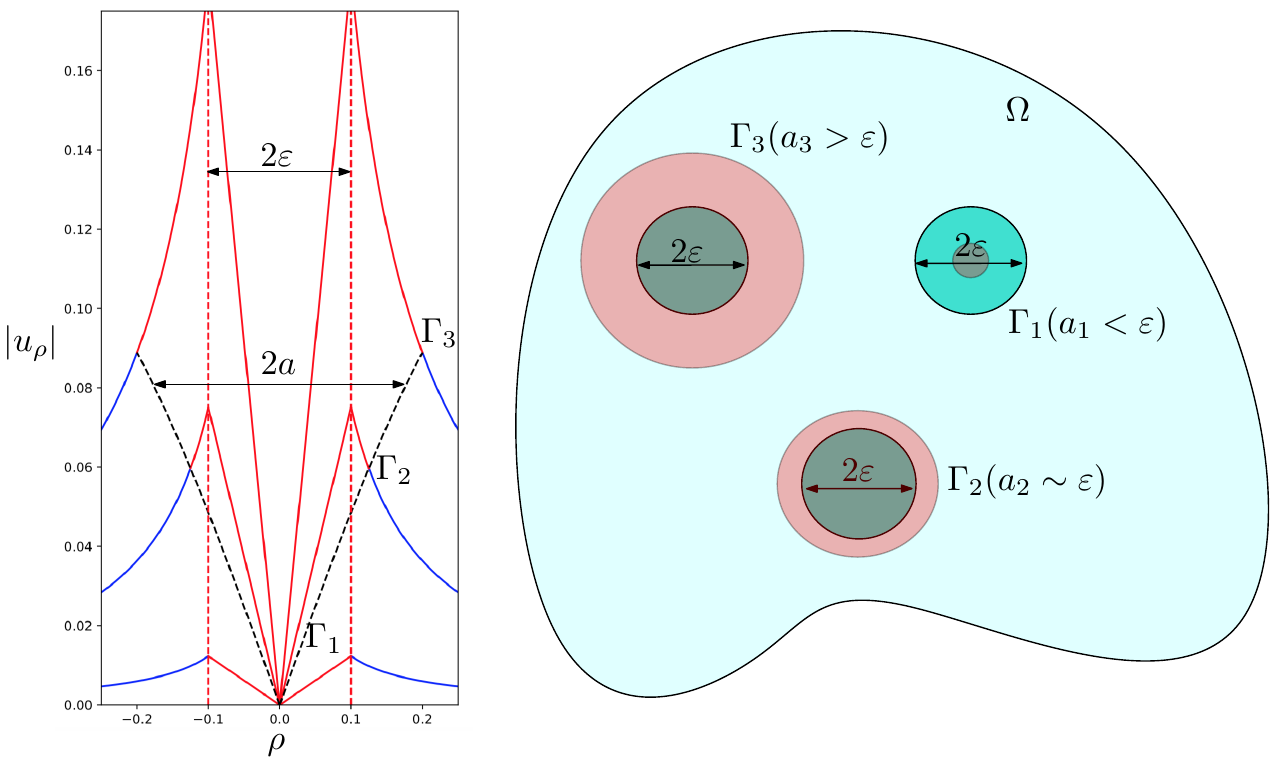}
  \caption{Example of vessel distribution with different radii $a_i$,
  and same fictitious area of influence $\varepsilon$.}\label{fig:approx-vessels}
\end{figure}

To tackle this issue, we generalize our approach one step
 further.  Namely, let us consider an additional parameter
 $\varepsilon>0$, representing the \emph{fictitious area of influence
   of the vessels}, and the corresponding circle $\Ve{\varepsilon}$ of
 radius $\varepsilon$. We now look for a fictitious elasticity problem
 inside both $\Ve{a}$ and $\Ve{\varepsilon}$, so that the solution
 coincides with the one defined in \eqref{eq:2d-exact} only outside of
 the ball of radius $\max\{ a,\varepsilon\}$ (see
 Figure~\ref{fig:approx-vessels}). Moreover, we impose a jump of the normal
 stress on the (non-physical) boundary $\Gamma_\varepsilon$ instead of
 on the vessel boundary $\Gamma_a$.

This fictitious problem can be constructed analogously to the one defined in Section \ref{ssec:2d-sing}.
Namely, defining  a continuous extension of the solution over the whole domain $\Omega$, computing the
jump of the normal stress across the boundary of $\Ve{\varepsilon}$ and imposing
this jump via a singular term in the elasticity equation.
In this case, the jump to be imposed reads
 \begin{equation}\label{eq:jump-epsilon}
      g_{\varepsilon}  := \frac{a^2}{\varepsilon^{2}} \frac{R^{2} \overline{p}
        \left(\lambda + 2 \mu\right)}{\left(R^{2} \mu + \lambda
          a^{2} + \mu a^{2}\right)} =
      \frac{a^2}{\varepsilon^{2}} g_a.
    \end{equation}

It is worth noticing that \eqref{eq:jump-epsilon} generalizes the formula previously derived in 
\eqref{eq:2d-sigma-jump}. In particular, the solution 
outside $\Ve{a}$ reduces to \eqref{eq:2d-extended} only if
$\varepsilon \leq a$. When $\varepsilon > a$, the solution coincides
with the exact one only outside of $\varepsilon$, and there is a
region, corresponding to the area between $\varepsilon$ and $a$, where
the solution is unphysical.

The advantage of using \eqref{eq:jump-epsilon}, is that  we can define arbitrarily the 
scale $\varepsilon$, that represents the \emph{resolution} of interest, i.e., the 
relevant scale at which we want to approximate our singular forcing term. 
At a distance at least $\varepsilon$ from the
vessels, the solution coincides with the expected one, while inside
the vessels, or inside a ball of radius $\varepsilon$ from the vessel
(whichever is bigger), the solution we obtain is unphysical.

The resolution at which we need to integrate over
$\Gamma_\varepsilon$ is now independent on the vessels size $a$, and,
in particular, it can be fixed \emph{a posteriori}, after a
discretization strategy (and a mesh size) is defined for the domain
$\Omega$.
This allows to define a forcing term in the limit for
$\varepsilon \to 0$, independently on the vessel's size. In this case,
the forcing term for a single vessel centered at the origin in the
variational formulation would reduce to
\begin{equation}
  \begin{aligned}
    \label{eq:hs-derivation}
    \lim_{\varepsilon\to 0} \int_{\Gamma_\varepsilon} g_{\varepsilon} \bn \cdot \vv \d \Gamma& =\lim_{\varepsilon\to 0}
    \int_{\Ve{\varepsilon}} \frac{a^2}{\varepsilon^{2}} g_a \ldiv\vv \d \bx
     \\
    & = \pi a^2 g_a \ldiv\vv
    (\boldsymbol{0}), \quad \forall \vv \in C^1(\Omega).
  \end{aligned}
\end{equation}

Equation~\eqref{eq:hs-derivation} defines the hyper-singular forcing
term
\begin{equation}
  \label{eq:hyper-singular-definition-delta}
   \bF^H(\bx)  :=  -\pi a^2 g_a \grad \delta(\bx),
\end{equation}
so that
\begin{equation}
  \label{eq:hyper-singular-definition}
  < \bF^H, \vv> :=  \pi a^2 g_a \ldiv\vv   (\boldsymbol{0}) \qquad \forall
  \vv \in C^1(\Omega).
\end{equation}

We remark here that $\bF^H$ cannot be used \emph{as-is} as a source
term for our elasticity problem, since it does not belong to the space
$H^{-1}(\Omega)$. 

It is however possible to \emph{mollify} the hyper singular
formulation \eqref{eq:hyper-singular-definition-delta}, by employing a
smooth approximation of the Dirac delta distribution
$\delta^{\varepsilon'}$, according to a small parameter
$\varepsilon'$, that again represents the resolution at which we
resolve our singular forcing terms.
Although this new parameter is technically
different from the one introduced in Equation~\eqref{eq:jump-epsilon},
in the rest of the paper we will set $\varepsilon' = \varepsilon$, i.e., 
identifying the scale of interest with the radius of approximation of the Dirac
delta distribution.

In particular, we consider approximations $\delta^\varepsilon$ of the
Dirac delta distribution such that:
\begin{itemize}
\item $\delta^\varepsilon(\bx-\by)  =  \delta^\varepsilon(\by-\bx)$
\item $\int_{\Re^2} \delta^\varepsilon(\bx-\by) \d \by =
  \int_{B_\varepsilon(\bx)} \delta^\varepsilon(\bx-\by) \d \by  = 1$
\item $\delta^\varepsilon \in C^1(\Re^2)$
\item $\int_{\Re^2} \nabla_\by \delta^\varepsilon(\bx-\by) \d \by =
  \int_{B_\varepsilon(\bx)} \nabla_\by \delta^\varepsilon(\bx-\by) \d \by = \bs{0}$.
\end{itemize}
and we defined the \emph{mollified} forcing term
\begin{equation}
  \label{eq:hyper-singular-mollified-definition}
  < \bF_\varepsilon^H, \vv> :=  \int_\Omega \delta^\varepsilon(\by)\pi
  a^2 g_a \ldiv\vv   (\by) \d \by \qquad \forall
  \vv \in H^1(\Omega).
\end{equation}
For a discussion on the properties of possible Dirac delta
approximations to use, we refer the reader to the excellent review
paper~\cite{Hosseini2014}.

The above formula can be straightforwardly generalized to the case of $N$ vessels, of radii
$a_i$, $i=1,\hdots,N$ and centered in $\bx_i$, $i=1,\hdots,N$. 
Introducing also the approximation
$a \ll R$ for the definition of the stress jump, i.e., 
\begin{equation}\label{eq:2d-sigma-jump-approx}
  g_a  \bn  = \frac{R^{2} \overline{p}
        \left(\lambda + 2 \mu\right)}{R^{2} \mu + \lambda a^{2}
        + \mu a^{2}} \bn  = \frac{2\mu + \lambda}{\mu} \overline{p} + O\left( \left( \frac{a}{R} \right)^2\right)\,,
\end{equation}
we obtain the hyper-singular forcing term
\[
\bF^H(\bx) = 
-\sum_{i=1}^N  \frac{2 \mu + \lambda}{\mu} \pi a_i^2 \overline{p}_i \, \grad \delta (\bx - \bx_i)\,.
\]
and its mollified version:
\[
\bF_\varepsilon^H(\bx) = 
-\sum_{i=1}^N \frac{2 \mu + \lambda}{\mu} \pi a_i^2 \overline{p}_i \, \grad \delta^\varepsilon (\bx - \bx_i)\,.
\]

In the two dimensional model, one can use either
$\bF^S_a, \bF^S_\varepsilon$, or $\bF^H_\varepsilon$ as forcing terms
and obtain a solution that approximates the exact solution outside
of the vessels up to higher order terms with respect to both the
ratio $a/R$ and $\varepsilon$.

When considering finite dimensional approximations, the first two
choices require the full discretization of the vessel boundary
$\Gamma$ or of the fictitious boundary
$\partial B_\varepsilon$, while employing $\bF^H_\varepsilon$
only requires evaluation of the integrals expressed in
Equation~\eqref{eq:hyper-singular-mollified-definition}.

\section{Three-dimensional case}
\label{sec:3d}
In three dimensions, we consider the vasculature as a network of
vessels, where each vessel is approximated as a thin
cylindrical domain, described via a one-dimensional manifold, denoted
as the \textit{centerline}, and a radius varying along the centerline.
In order to obtain the singular source terms, we will then integrate
over the centerline the equivalent of the two-dimensional formulation
discussed in Section \ref{sec:2d} considered on a plane that is
locally orthogonal to the centerline.

In what follows, we will also assume that curvature of each vessels, within a single segment,
varies slowly w.r.t. to its arclength, so that its effect, as well as
elastic effects of the vessels, may be neglected. For a possible way
to include the elastic behaviour of the vessels we refer
to~\cite{alzettaheltai-2018-a}. In each cross sectional plane of the
vessel, we approximate the local behaviour of the problem as in the
two-dimensional axi-symmetric case.

\subsection{Geometrical setting}
\begin{figure}
\centering
\includegraphics{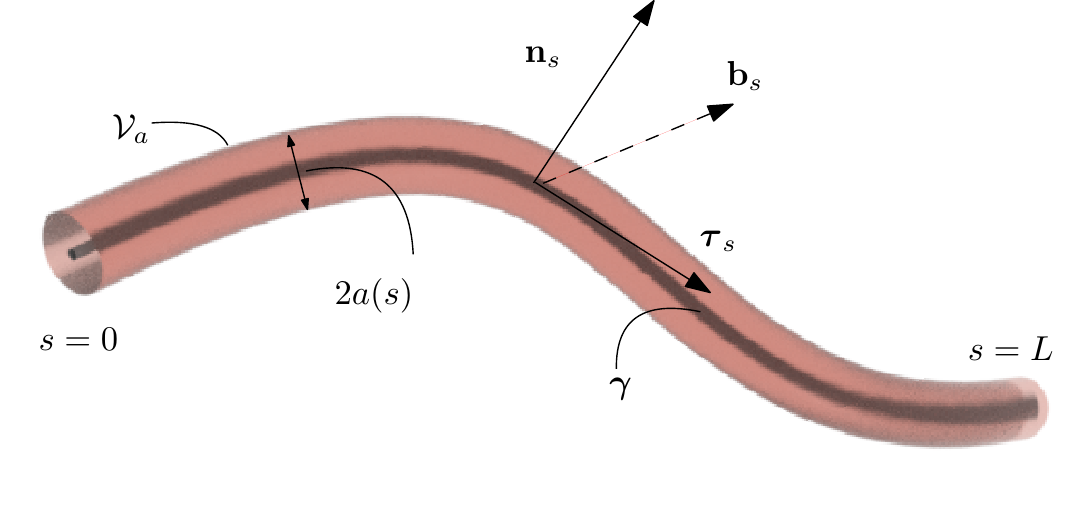}
\caption{Sketch of a 3D vessel, identified via its cross-sectional radius and its  centerline.}
\label{fig:vessel-3d}
\end{figure}

In order to introduce the geometrical model of vascular network, we decompose the network in a set of
non-intersecting vessel segments. For each segment, let us introduce a one-dimensional arc-length curve
\[
\bgamma (s): [0,L] \to \Omega \subset R^3,
\]
describing the vessel centerline, and a positive function
\[
a(s): [0,L] \to \Omega \subset R,
\]
standing for the radius of the cross-section at each $s \in [0,L]$.
Moreover, let us denote with $A(s)$ the cross-section, i.e., 
the disk of radius $a(s)$ orthogonal to $\bgamma(s)$, and with $|A(s)| = \pi a^2(s)$
the cross-sectional area, for all $s \in [0,L]$
(a sketch depicting these quantities is provided in Figure \ref{fig:vessel-3d}).

In order to formally derive the multiscale model, we introduce at each $s \in [0,L]$ also the
\textit{Frenet} frame $\btau_s = \gamma'(s)$ (tangential vector in $s$), 
$\mathbf n_s=\btau_s'/|\btau_s'|, \mathbf b_s= \btau_s\times \mathbf n_s$ 
(basis of the normal plane in $s$).

Finally, we assume that the fluid pressure within the vessel is constant over each cross section, thus
introducing a function
\[
p_{\bgamma}(s) : [0,L] \to \Omega \subset R,
\]
denoting the excess pressure along the centerline.

Let us consider a three-dimensional domain $\Omega$, the set
\begin{equation}\label{Va-def}
\vessel_a(\bgamma) = \{ \bx \in \Omega \text{ s.t. } \text{dist}(\bx,
\bgamma) < a \}, 
\end{equation}
(denoting the vessel domain) and the tissue domain
$\Omega_a = \Omega \backslash \vessel_a$. 
Moreover, let $\Gamma = \partial \vessel_a(\bgamma)$.

Assuming that the domain $\vessel_a(\bgamma)$ describes a non-intersecting vessel segment, the coordinate transformation mapping
\begin{equation}
  \label{eq:coordinate-transformation}
  \bphi(r, \theta, s) := \bgamma(s) + r\cos(\theta) \mathbf n_s +
  r\sin(\theta) \mathbf b_s,
\end{equation}
is one-to-one from a cylindrical domain in polar coordinates
$(r,\theta,s) \in (0,a(s)] \times [0,2\pi] \times [0,L]$ onto
$\vessel_a(\bgamma)$.

We denote with $\gamma^{-1}: \vessel_a(\bgamma) \mapsto [0,L]$ the
function that identifies, for each point $\bx$ in
$\vessel_a(\bgamma)$, the arc-length coordinate $s \in [0,L]$ such that
$\bgamma(s)$ has minimum distance from $\bx$, i.e.,
\begin{equation}
  \label{eq:definition-bgamma-inverse}
  \gamma^{-1}(\bphi(r,\theta,s)) := s, \qquad \forall r \in [0,a(s)],
  \quad \forall
  \theta \in [0,2\pi].
\end{equation}

Assuming that the pressure is constant over cross-sections, we define
an extension $p$ (defined in the three dimensional vessel) of the 
one-dimensional excess pressure $p_{\bgamma}$ via
\[
p(\bx) = p_{\bgamma}\left(\gamma^{-1}(\bx)\right)\,, \forall \bx \in \vessel_a\,.
\]

\subsection{The singular formulation}

As in the previous case, we aim to solve an elasticity problem on the domain $\Omega_a$
by constructing variational
formulation on the whole three-dimensional domain $\Omega$, in which the elasticity
problem in $\Omega_a$ is extended by a fictitious problem in $\vessel_a$. 
Proceeding analogously as in Section \ref{ssec:2d-sing}, we seek for the solution of 
a problem of the form

\begin{problem}[3D, variational]
  \label{pb:three-dimensional-variational}
  Given the excess pressure field $p_{\bgamma}$, the vessel
  configuration $\bgamma$, and the radius function $a$, find
  $\bu \in V$ such that
 \begin{equation*}\label{eq:force-3d}
    (2\mu \tens{e}(\bu), \tens{e}(\vv) )_{\Omega} + (\lambda \ldiv\bu, \ldiv\vv)_{\Omega}  = 
    <\bF_{(\bgamma, p_{\bgamma}, a)},\vv> \quad
    \forall \vv \in V.
  \end{equation*}
\end{problem}

The source term $\bF_{(\bgamma, p_{\bgamma},a)}$ shall be defined in such a way to enforce, for
each $s \in [0,L]$, a given jump of the normal stresses across
$\Gamma \cap A(s)$.

Ideally, we would like to use 
the same reasoning that lead to the definition of the two-dimensional model
problem, that is, defining
$\bF_{(\bgamma, p_{\bgamma}, a)}$ such that the solution $\bu$ coincides
with the one that would be obtained by solving the true problem in the
elastic matrix alone, with non-homogemous Neumann boundary conditions
on the vessel boundary $\Gamma$ depending on the pressure $p_{\bgamma}$.

However, in three dimensions, an explicit solution is only available
for trivial vessel geometries and boundary conditions, unless we
assume that all quantities that change along the vessel coordinate
direction varies slowly w.r.t. to $s$. In this case we could still use
the same principle used in the two dimensional approximation by
integrating the derivation of the two-dimensional model problem along
the arclength $s$.
We start by constructing a force distribution $\bF^S(\bgamma, p_{\bgamma}, a)$ given by
\begin{equation}
  \begin{aligned}
    <\bF^S_{(\bgamma, p_{\bgamma}, a)},\vv> & := \int_{\Gamma} \frac{(2\mu+\lambda)}{\mu} p(\by) \, \bn(\by)\cdot \vv(\by) \d \Gamma_{\by}
    \\ & =  \int_{\Gamma} \frac{(2\mu+\lambda)}{\mu} p(\by) \, \bn(\by)\cdot \int_{\Omega} \vv(\bx) \delta(\bx-\by) \d \bx
    \d \Gamma_{\by} \\
    & =  \int_{\Omega} \vv(\bx) \cdot \left(
      \int_{\Gamma} \frac{(2\mu+\lambda)}{\mu} p(\by) \, \bn(\by)   \delta(\bx-\by) 
      \d \Gamma_{\by}\right)  \d \bx\,.
  \end{aligned}
\end{equation}
In 3D, it is therefore possible to define the singular
source term as an integral over the vessel boundary of the
form
\begin{equation*}
              \bF^S_{(\bgamma, p_{\bgamma},a)}(\bx) := \int_{\Gamma}
          \frac{(2\mu+\lambda)}{\mu} p(\by)\, \bn(\mathbf y) \delta(\bx-\mathbf y) \d
          \Gamma_{\mathbf y},\; \forall \bx \in \Omega\,.
  \end{equation*}

The next step is to generalize the singular formulation in order to impose a given jump of the normal stresses
across the boundary of a vessel defined by the centerline $\bgamma([0,L])$ and an
arbitrary, constant, radius $\varepsilon > 0$.
Let $\vessel_{\varepsilon}(\bgamma)$ denote the generalized vessel of
radius $\varepsilon$, defined analogously to $\vessel_a(\bgamma)$ in
\eqref{Va-def}, and let
$\Gamma^\varepsilon = \partial V_{\varepsilon}(\bgamma)$ represent the
boundary of such generalized vessel.

Proceeding as in the two-dimensional case, we introduce 
\begin{equation}
  \label{eq:definition-of-ga}
  \hat g_a(s) := \pi a^2(s) p_{\bgamma}(s) \frac{(2\mu+\lambda)}{\mu}.
\end{equation}
which is based on the approximation of the jump across the vessel
boundary in the two-dimensional case (see
\eqref{eq:2d-sigma-jump-approx}), and the forcing term
\begin{equation}
  \label{eq:Fepsv3d}
  \bF^S_{(\bgamma, p_{\bgamma},\varepsilon)} (\bx) := \int_{\Gamma^\varepsilon}
  \frac{\hat  g_a(\gamma^{-1}(\by))}{\pi\varepsilon^2}\, \bn(\mathbf y) \delta(\bx-\mathbf y) \d
  \Gamma_{\mathbf y},\; \forall \bx \in \Omega\,.
\end{equation}

\subsection{The hyper-singular formulation}

We introduce the gradient operator in the plane orthogonal to $\btau$ as
\begin{equation}\label{eq:grad-tau}
\nabla_{\btau} \bu := \left( \tens{1} - \btau \otimes \btau \right) \nabla \bu
\end{equation}
and the planar divergences as 
\begin{equation}\label{eq:div-tau}
\nabla_{\btau} \cdot \bu := \text{tr} \left( \nabla_{\btau} \bu \right)=
\nabla \cdot \bu  - \btau \cdot \left( \nabla \bu \, \btau\right)\,.
\end{equation}

\begin{remark}
  Notice that, in the case $\btau= (0,0,1)$ (vessel directed orthogonal to the
  $(x,y)$-plane), the above definition
  reduces to the gradient and divergence operators considered for the
  two-dimensional case.
\end{remark}

For easiness of notation, let us denote
\begin{equation}
  \label{eq:extension-ga}
   \hat  g := \hat  g_a \circ \gamma^{-1}, 
\end{equation}
i.e., the extension of $\hat  g_a$ on $\vessel_\varepsilon(\bgamma)$.

For any function $\vv \in (C^1(\Omega))^3$ the singular force can be rewritten as
\begin{equation}\label{eq:Feps3d-0-0}
\begin{aligned}
 <\bF^S_{(\bgamma, p_{\bgamma},\varepsilon)},\vv> & = 
 \int_{\Gamma^\varepsilon} \frac{\hat  g}{\pi
    \varepsilon^2} \vv\cdot \bn \d \Gamma_{\bx} =   
    \int_{\partial\vessel_\varepsilon(\bgamma)} \frac{\hat  g}{\pi
    \varepsilon^2} \vv\cdot \bn \d \Gamma_{\bx} - \int_{A_{\varepsilon}(0)} \frac{\hat  g}{\pi
    \varepsilon^2} \vv\cdot \bn - \int_{A_{\varepsilon}(1)} \frac{\hat  g}{\pi
    \varepsilon^2} \vv\cdot \bn
    \\
    & = \int_{\vessel_\varepsilon(\bgamma)}
  \frac{1}{\pi \varepsilon^2} \nabla \cdot \(\hat g \vv\) \d
  \bx + \int_{A_{\varepsilon}(0)} \frac{\hat  g}{\pi
    \varepsilon^2} \vv\cdot \btau_0 - \int_{A_{\varepsilon}(1)} \frac{\hat  g}{\pi
    \varepsilon^2} \vv\cdot \btau_1
      \end{aligned}
\end{equation}
denoting with $A_{\varepsilon}(0)$ and $A_{\varepsilon}(L)$ the bottom
and the top face of the vessel, respectively, and noticing that
$\btau_0 = - \bn$ on the bottom face (as $\btau$ is directed along the
vessel, while $\bn$ is directed outwards).

We now consider the limit of $\bF^S_{(\bgamma, p_{\bgamma},\varepsilon)}$ for
$\varepsilon\to 0$, approximating $\vv$ within the vessel with $\vv \circ \gamma^{-1}$, i.e., with its 
value on the vessel centerline. In view of \eqref{eq:Feps3d-0-0} we obtain
\begin{equation}\label{eq:Feps3d-0-1}
\begin{aligned}
   \lim_{\varepsilon \to 0} <\bF^S_{(\bgamma, p_{\bgamma},\varepsilon)},\vv> & = 
   \lim_{\varepsilon \to 0} \int_0^L\int_0^{2\pi}\int_0^\varepsilon
  \frac{1}{\pi \varepsilon^2} \nabla \cdot \(\hat  g \vv\) \d r \, r\d\theta
  \d s  + \underbrace{ \hat g \vv \btau_0 - \hat g \vv \btau_L}_{=- \int_0^L \frac{\partial}{\partial s} \hat g \vv \btau} \\
  & = \lim_{\varepsilon \to 0} \int_0^L \nabla \cdot \(\hat  g \vv\) \d s - 
  \int_0^L \btau \cdot \nabla ( \hat g \vv \circ \bgamma ) \btau \d s\,.
   \end{aligned}
\end{equation}
Hence, using the definition \eqref{eq:div-tau} and observing that 
\[
\int_0^L \frac{\partial}{\partial s} \hat g \vv \btau = \int_0^L \btau \cdot \nabla ( \hat g \vv \circ \bgamma ) \btau \d s
\]
yields
\begin{equation}\label{eq:Feps3d-0}
\begin{aligned}
   \lim_{\varepsilon \to 0} <\bF^S_{(\bgamma, p_{\bgamma},\varepsilon)},\vv> & = 
    \lim_{\varepsilon \to 0}  \int_{\vessel_\varepsilon(\bgamma)} \frac{1}{\pi \varepsilon^2} \nabla_{\btau}\(\hat g \vv\) \d   \bx 
  + \int_0^L \btau \cdot \nabla ( \hat g \vv \circ \bgamma ) \btau \d s 
- \int_0^L \btau \cdot \nabla ( \hat g \vv \circ \bgamma ) \btau \d s \\
& = \int_0^L \nabla_{\btau}\(\hat g \vv\) \d  s
= \int_{\Omega} \left[\int_0^L  \nabla_{\btau}\(\hat g \vv \,\delta(\bx - \bgamma(s)\) \d s\right]\d \bx,
\qquad \forall \vv \in (C^1(\Omega))^3.
\end{aligned}
\end{equation}

%
In view of \eqref{eq:Feps3d-0}, we consider the variational formulation 
 \begin{equation*}
    (2\mu \tens{e}(\bu), \tens{e}(\vv))_{\Omega} + (\lambda \ldiv\bu, \ldiv\vv)_{\Omega} = 
    <\bF_{(\bgamma, p_{\bgamma},a)}^{H} + \bF_{(\bgamma, p_{\bgamma},a)}^{\btau} ,\vv>
    \quad
    \forall \vv \in C^1(\Omega)
  \end{equation*}
where the right hand side can be defined through the hyper-singular term
              \begin{equation}\label{eq:3d-hyper}
    \bF_{(\bgamma, p_{\bgamma},a)}^{H}(\bx)  := \int_0^L \hat g_a(s)
     \grad_{\btau} \delta(\bx-\bgamma(s)) \d s,\; \qquad \forall \bx \in \Omega\,
  \end{equation}
and the singular source
\begin{equation}\label{eq:3d-tgt}
\bF_{(\bgamma, p_{\bgamma},a)}^{\btau}(\bx)=
\int_0^L  \hat g_a'(s)  \, \delta(\bx-\bgamma(s)) \btau \d s
,\; \qquad \forall \bx \in \Omega\,.
\end{equation}
 which has support on the centerline and it is directed tangential to it.
 In particular, if vessel radius and pressure are constant along $\gamma$, the 
 singular term \eqref{eq:3d-tgt} vanishes, and the immerser method reduces to a
 hypersingular force equal to the tangential derivative of a Dirac delta function.
 
%

  %

\begin{remark}
Notice that the forces introduced in \eqref{eq:3d-hyper}--\eqref{eq:3d-tgt}
depend only on one-dimensional information, such as
centerline, the excess pressure $p(s)$, the radius, and the cross-sectional area, and it allows therefore to
represent the vessel uniquely through a one-dimensional manifold.
\end{remark}

\begin{remark}
  Similarly to what happens in the two-dimensional case, the
  forcing terms are
  not in $V^*$, and we should replace $\bF_{\bgamma}^H$ and $\bF_{\bgamma}^{\btau}$
  with a mollified version, where the
  Dirac delta distribution $\delta$ is replaced by a regularized
  version of it,  depending on a (small) regularization parameter (see Section \ref{ssec:discrete}).
\end{remark}

\section{Homogenized behavior of pressurized tissues}\label{ssec:tissue-homo}
Let us consider a tissue domain $\Omega$, with a given one-dimensional characterization
of the vasculature (centerlines, radii, and pressures), so that
the immersed method based on the singular forces \eqref{eq:3d-hyper}-\eqref{eq:3d-tgt} can be defined.
In pactice, this information can be either (fully or partially) recovered  from medical imaging 
(e.g., diffusion MRI) or generated artificially using statistical methods (as it will be shown later
in Section \ref{ssec:3dreal}). Let $\beta$ represents the volume fraction of $\Omega$ that is covered by
vessels (e.g., for soft tissues, $\beta$ is typically below 5\%).
Using the structure of the immersed finite element method, 
we can additively decompose the solution, isolating the effect of the
pressurized vessels on the right hand side. In other words, 
we can seek the solution $\bu^p_h$ in $V_h$ such that
\begin{equation}
  \label{eq:homogenous-vascularized}
  (\tens{\sigma}(\bu^h_p), \tens{e}(\vv_h)) = (\bv{F}(p, a, \beta), \vv_h)_\Omega \qquad
  \forall \vv_h \in V_h.
\end{equation}
Starting from Equation \eqref{eq:homogenous-vascularized}, 
the goal of this section is to derived a homogenized characterization of the effective mechanical 
properties of a pressurized tissue, depending on the properties (geometry, density, pressure)
of the underlying vasculature.

\subsection{Derivation in the two-dimensional case}\label{ssec:derivation}
To begin with, let us consider the two-dimensional problem, 
which can be seen as a cross section of
a three-dimensional domain, in the case that all vessels are directed along the $z$-direction.
We assume an
uniform random spatial distribution of $n$ random vessels located in 
$\{\bx_i\}_{i=1}^n$, with radius $a_i$ and with excess
pressure $p_i$, for $i = 1,\hdots,n$.

\begin{figure}
  \centering
  \includegraphics{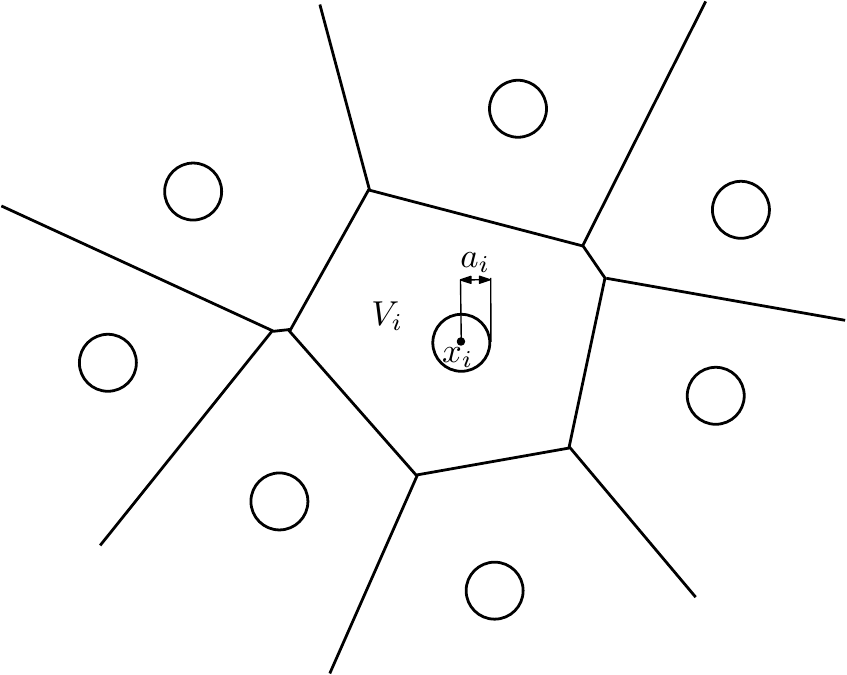}
  \caption{Voronoi diagram of the vessel centers for the
    two-dimensional case.}
  \label{fig:voronoi}
\end{figure}

Let us now consider the Voronoi diagram covering the domain
$\Omega$ with generators in $\{\bx_i\}_{i=1}^n$ (as in
Figure~\ref{fig:voronoi}). Indicating with $V_i$ the $i$-th Voronoi
cell, and with $|V_i|$ its volume, we can define the local vessel
density $\beta_i = \pi^2a_i/|V_i|$, and we can interpret the forcing
term in equation \eqref{eq:homogenous-vascularized} as the
approximation through the Voronoi diagram of a \emph{continuous}
integral over the domain $\Omega$:
\begin{equation}
  \label{eq:simplified-f-uniform}
  \begin{aligned}
  <\bv{F},\vv> & := \sum_{i=1}^{n}\frac{(2\mu+\lambda)}{\mu}p_i 
  \pi a_i^2 \ldiv\vv(\bx_i) = \sum_{i=1}^{n} |V_i| \frac{(2\mu+\lambda)}{\mu}p_i
  \frac{\pi a_i^2}{|V_i|} \ldiv\vv(\bx_i)  \\
  & = \sum_{i=1}^{n}  \beta_i |V_i|  \frac{(2\mu+\lambda)}{\mu}p_i \ldiv\vv(\bx_i) 
  \\ & \simeq \sum_{i=1}^{n} \int_{V_i} \beta_i \frac{(2\mu+\lambda)}{\mu}p_i \ldiv\vv(\bx)
  \simeq \int_\Omega \beta \frac{(2\mu+\lambda)}{\mu}
  p \ldiv\vv \d \Omega,
    \end{aligned}
\end{equation}
where $\beta$ and $p$ are homogenised quantities that can vary
spatially, and represent the local excess pressure and the local
vessel density of the tissue.

This approximation yields an homogenized elasticity problem
\begin{equation}
  \label{eq:homogenous-vascularized-2}
  (\tens{\sigma}(\bu^h_p), \tens{e}(\vv_h)) = \int_\Omega \beta \frac{(2\mu+\lambda)}{\mu}
  p \ldiv\vv_h \d \Omega  \qquad
  \forall \vv_h \in V_h,
\end{equation}
where the forcing term on the right-hand-side acts on the dilatational part of the
deformation. The solution to \eqref{eq:homogenous-vascularized-2}
satisfies the following conservation equation:
\begin{equation}
  \label{eq:conservation}
  \int_\Omega  2\mu |\tens{e}(\bu)|^2 + \lambda |\ldiv\bu|^2 \d \Omega = \int_\Omega  \beta \frac{(2\mu+\lambda)}{\mu}
  p \ldiv\bu \d \Omega,
\end{equation}
where one clearly sees that the pressurized vessel network
acts as a non-conservative pressure source in the energy conservation
equation.

Let us now assume homogeneous Neumann boundary conditions on the
outer boundary. The force due to the uniform distribution of vessels (with
constant pressure and fixed volume fraction) 
would produce -- up to rigid deformations -- a uniform
dilation (or compression). We seek therefore a solution of the elasticity problem of the form
\begin{equation}
  \label{eq:exact-solution-uniform-pressure}
  \bu = c \bx\,.
\end{equation}

In order to determine the constant $c$, we insert
\[
|\ldiv\bu|^2 = (2c)^2  = 4 c^2,\;
|\tens{e}(\bu)|^2 = 2c^2 
\]
into \eqref{eq:conservation}, obtaining
\[
4\mu c^2  + 4\lambda c^2  = 2 \beta \frac{(2\mu+\lambda)}{\mu}p c
\]
and thus
\begin{equation}
  \label{eq:c-val}
 c =  \frac{\beta p}{2\mu}\, \frac{2\mu+\lambda}{\mu  + \lambda}.
\end{equation}

 The effect of the pressurized vasculature produces therefore a stress on the boundary equal to
\begin{equation}
  \label{eq:expected-pressure-vascularized}
  \tens{\sigma}(\bu^h_p) \bn = \frac{\beta p}{\mu} (2\mu+\lambda)\bn,
\end{equation}
and the corresponding total force $\bF_p$ on the face $A$ can be
computed as
\begin{equation}
  \label{eq:estimated-total-pressure-force}
  \bF_p :=\int_A \tens{\sigma}(\bu^h_p) \bn  \d A = |A|\frac{\beta p}{\mu} (2\mu+\lambda)\bn.
\end{equation}

In the case of a uniform spatial distribution of vessels, with constant vessel sizes
and constant pressure, we obtain $(\bF_p \cdot \bn) \bn = \bF$, i.e.,
the internal force generated by the pressurized vasculature is always
directed along the normal direction to the surface.

Notice that the derivation of the total force
\eqref{eq:estimated-total-pressure-force} is based on the assumptions
of  Neumann boundary conditions and on the fact
that the term $\beta\frac{(2\mu+\lambda)}{\mu}p$ is constant across
the domain. In this situation, the divergence theorem yields
\begin{equation}
  \label{eq:simplified-f-uniform-byparts}
  <\bv{F},\vv> = \int_\Omega \beta \frac{(2\mu+\lambda)}{\mu}
  p \ldiv\vv \d \Omega = \int_{\partial\Omega} \beta \frac{(2\mu+\lambda)}{\mu}
  p \vv \cdot \bn \d \Gamma.
\end{equation}
The same argument cannot be used in the case of homogeneous Dirichlet boundary conditions.
In fact, the term on the right-hand-side of
Equation~\eqref{eq:simplified-f-uniform-byparts} would be tested
against functions $\vv$ in $V$, whose value on $\partial \Omega$ would
be zero, meaning that a uniform distribution
of vessels with constant pressure has no effect on the
solution.

In reality, since the distribution of vessels is discrete (although uniform), its
effect should be noticeable also with Dirichlet boundary conditions,
by measuring $\bF_p = \int_A \tens{\sigma}(\bu^h_p) \bn$ for each of the
faces of the domain. By linearity, this should be equal (on average)
to Equation~\eqref{eq:estimated-total-pressure-force}.

\subsection{Derivation in the three-dimensional case}\label{ssec:derivation-3d}
In order to generalize the above arguments to the three-dimensional
case, let us first consider a uniformly pressurized tissue where the
distribution of vessels is spatially uniform and only aligned in a
specified direction $\btau$. In this case, the force exerted by the
presence of the vessels is isotropic in the plane which is orthogonal
to the vessels direction.

The singular force due to the presence of a constant pressure and
uniform distribution of vessels, can then be written as
\[
<\bv{F},\vv> = \int_\Omega p \frac{(2\mu+\lambda)}{\mu}
 \beta \nabla_{\btau} \cdot \vv \d \Omega = \int_{\partial \Omega} p \frac{(2\mu+\lambda)}{\mu} \beta\left( \tens{1}- \btau \otimes \btau \right) \vv \cdot \bn\,,
\]
resulting in a uniform internal stress in the material, given by
\begin{equation}\label{eq:sigma-tau}
  \sigma(\bu)  = p \frac{(2\mu+\lambda)}{\mu} \beta \left(\tens{1} - \btau \otimes \btau \right).
\end{equation}

To investigate the effect of the internal pressure stress, we observe
that, by symmetry considerations, a solution with homogeneous Neumann
boundary conditions should have the form
\begin{equation}\label{eq:u3d}
\bu = c \bx + d \btau \otimes \tau \bx
\end{equation}
(i.e., an uniform dilation/compression plus a uni-axial deformation
along the direction of the vessels).

In particular it holds:
\[
\tens{e}(\bu) = c \tens{1} + d \btau \otimes \btau, \qquad \nabla \cdot \bu = 3c + d
\]
which yields
\[
\sigma(\bu) = 2 \mu c  \tens{1} + 2 \mu \btau \otimes \btau + \lambda(3c + d)  \tens{1}.
\]

Equation~\eqref{eq:sigma-tau} indicates that the stress along the direction
$\btau$ vanishes, i.e., $\sigma(\bu) \btau = \boldsymbol{0}$, yielding
\begin{equation}\label{eq:d_of_c}
2\mu c + 2 \mu d +  3\lambda c + \lambda d = 0 \Rightarrow
 d = -\frac{2 \mu + 3 \lambda}{2\mu + \lambda}c
 \end{equation}

 Inserting \eqref{eq:d_of_c} into \eqref{eq:u3d} we obtain
 \[
 \bu = c \left(  \tens{1}- \frac{2 \mu + 3 \lambda}{2\mu + \lambda} \right)\bx
 \]
and 
\begin{equation*}
\begin{aligned}
\sigma (\bu) & = 2\mu \tens{e}(\bu) +\lambda  \nabla \cdot \bu \tens{1}=
c \left( 2\mu \left[  \tens{1} - \frac{2 \mu + 3 \lambda}{2\mu + \lambda} \btau \otimes \btau\right]
+ \frac{4\mu}{2\mu + \lambda} \tens{1}
\right) \\
& = c \frac{2\mu}{2\mu + \lambda} \left( 2\mu + 3 \lambda\right)( \tens{1}- \btau \otimes \btau).
\end{aligned}
\end{equation*}

The latter equation, combined with \eqref{eq:sigma-tau} when $\bn \perp \boldsymbol \tau$, yields 
$c = \frac{\beta p}{2\mu^2} \frac{(2\mu + \lambda)^2}{2\mu + 3 \lambda}$, which combined with
\eqref{eq:d_of_c} leads to the following solution for the displacement
\begin{equation}\label{eq:u-ex-3d}
\bu = \frac{\beta p}{2 \mu^2}(2\mu + \lambda)\left( \left( \frac{2\mu + \lambda}{2\mu + 3 \lambda} \right) \tens{1}- \btau \otimes \btau\right)\,.
\end{equation}
The latter,  together with~\eqref{eq:sigma-tau}, characterizes pressurized
tissue materials with a single uniform distribution of vessels with
volume fraction $\beta$, directed along the direction $\btau$, and
pressurized with constant pressure $p$.

\begin{figure}
  \centering
  \includegraphics[width=.9\textwidth]{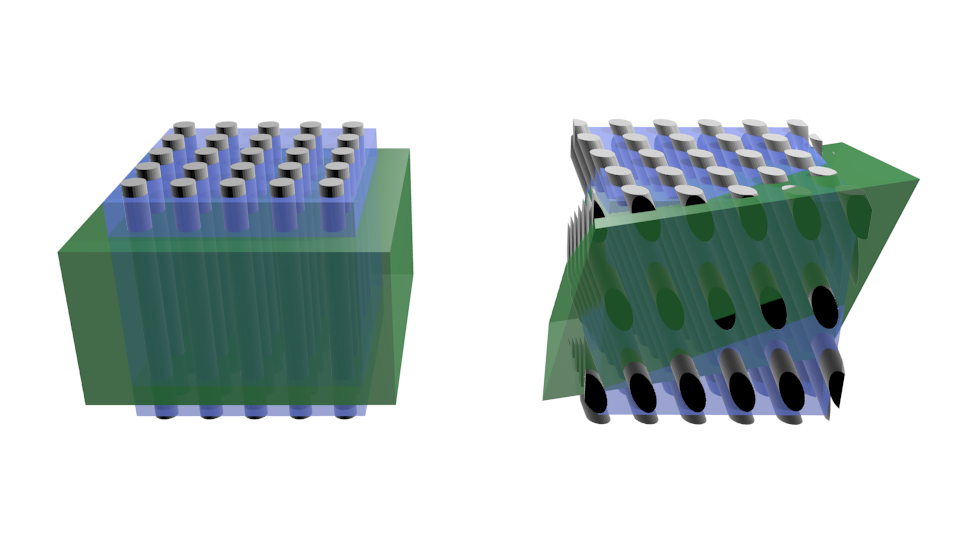}
  \caption{The effect of vessels orientation in a uniformly
    vascularized material (represented in blue in the figure), with a
    single dominant vessel direction: when vessels are aligned with
    one of the sample principal directions (left), the resulting
    deformation (represented in green in the figure) is orthogonal to
    the sample boundaries. If the vessels distribution is not aligned
    with one of the sample principal directions (right), a pressure
    induced shear is observed (right).}
  \label{fig:pressure-induced-shear}
\end{figure}

The material responds to pressurisation with an anisotropic
deformation: a contraction along the direction of the vessels, and a
dilation in the direction orthogonal to the vessels. When the
direction $\btau$ is not parallel to one of the axes, the anisotropic
contractile behaviour induces shear in the material, as depicted in
Figure~\ref{fig:pressure-induced-shear}.

Thanks to linearity, the model generalizes easily to uniform
distributions of vessels with more than one dominant direction. In
particular, an anisotropic vessel distribution can be described by a
tensor $\tens{\beta}$
 \begin{equation}
   \label{eq:anisotropic-distribution-density}
   \tens{\beta} := \sum_{i=1}^3 \beta_i \left(\tens{1} - \btau_{s_i} \otimes \btau_{s_i} \right),
 \end{equation}
 where the directions $\btau_{s_i}$ are mutually orthonormal, and
 $\beta_i$ represent the (constant) vessel spatial densities across
 planes that are orthogonal to the directions $\btau_{s_i}$.

The corresponding source term takes the form
\begin{equation}
  \label{eq:anisotropic-beta-source}
  <\bv{F}_\beta,\vv> = \int_\Omega p \frac{(2\mu+\lambda)}{\mu} 
  \text{tr} (\tens{\beta} \nabla \vv)\d \Omega = \int_{\partial
    \Omega}  p \frac{(2\mu+\lambda)}{\mu} \vv \cdot \tens{\beta}\bn,
\end{equation}

And the generated cauchy stress due to the pressurisation is given by
\begin{equation}
  \label{eq:cauchy-anisotropic-uniform-distribution}
  \sigma(\bu)  = p \frac{(2\mu+\lambda)}{\mu} \tens{\beta}.
\end{equation}

The symmetric tensor $\tens{\beta}$ may be used to characterize the
influence of pressurized vessels inside the tissue. In particular, for
densely distributed networks of vessels, one may approximate the
tensor $\tens{\beta}$ by the integral
\begin{equation}
  \label{eq:measuring-beta}
  \tens{\beta} := \frac{1}{|\Omega|} \int_\Gamma \pi a^2 (\tens{1} -
  \btau\otimes\btau)\d \Gamma,
\end{equation}
where $\Gamma$ is the centerline of the vessel network (the union of
all vessel segments), $\btau$ is the local tangent vector, and $a$ is
the local radius of the vessel. When the vessel network is not
explicitly available, because data resolution does not allow to
reconstruct it, one could infer the average properties of the
pressurized material (Lam\`e parameters, principal directions, and
local volume fraction), by a sequence of pure shear and pure dilation
measurements.

\subsection{Characterization of vascularized tissue properties}\label{ssec:tissue-charact}
In the context of magnetic resonance elastography, tissue characterization is based on
the solution of an inverse elasticity problem where the spatial resolution of available data
(i.e., displacement field acquired via phase contrast MRI) is of the order
of millimeters, and it is typically much coarser than
the scale of vascular structures (vessel diameters).
Neglecting the effect of vascular pressure might drastically change the values of the estimated parameters.
This was shown experimentally, e.g., in \cite{Chatelin2011}, where the shear modulus values
obtained via elastography ex-vivo were much lower than those found in vivo.

These observations demonstrate that the inverse modeling of tissue should be based
on effective material models (at the scale of available data) that are able to 
capture the influence of microscopic vasculature (and related pressures) on the
coarse mechanical parameters (Lam\'e coefficients).
%
This section discusses the implication of the homogenized model including the singular forces
(derived in Section \ref{ssec:tissue-homo}) in the characterization of mechanical properties of tissues.

Typical experimental settings -- targeted to the quantification of
elastic parameters -- are designed to induce one of two ideal
deformations: \emph{pure shear}, used to obtain information about the
shear modulus $\mu$, and \emph{free compression}, used to obtain
information about the Poisson ratio $\nu$, i.e., the ratio of relative
contraction to relative expansion of the material. The second Lam\'e
coefficient can then be extracted by the relation
$\lambda = \frac{2\mu\nu}{1-2\nu}$.

Pure shear experiments are easier to reproduce in \emph{in-vivo}
tissues, and mimic an essentially two-dimensional configuration where
the displacement is given by
$\bu_{ij} = \frac{c}2 (\bv{e}_i\otimes\bv{e}_j) \bv{x}$, where $c$ is
a controlled constant, given by the experimental setting, and
$i \neq j$. In these cases, it is easy to show that the corresponding
stress is given by
\begin{equation}
  \label{eq:pure-shear}
  \tens{\sigma}(\bu_{ij}) = c \mu_{ij} (\bv{e}_i\otimes\bv{e}_j + \bv{e}_j\otimes\bv{e}_i). 
\end{equation}

A measure of the force along the $i$-axis, measured on the face with
surface $|A|$ and normal $\bv{e}_j$ gives:
\begin{equation}
  \label{eq:pure-forces}
  \begin{split}
    F_{ij} & := \int_{A_j} (\tens{\sigma}(\bu_s) \bv{e}_j)\cdot
    \bv{e}_i \d\partial \Omega = c |A| \mu_{ij},
  \end{split}
\end{equation}
that provides direct access to the shear modulus
$\mu_{ij} = F_{ij}/(c|A|)$. In case of isotropic materials,
$\mu_{ij} = \mu$ on every face, and one experiment is enough to
characterize the elastic matrix.

However, when the tissue contains pressurized fluid vessels
characterized by the density distribution $\tens{\beta}$, the internal
stress does not contain only the terms in
equation~\ref{eq:pure-shear}. In this situation, the homogeneized
characterization introduced in Section~\ref{ssec:tissue-homo}
shows that the pressure induced stress is in general
anisotropic, and adds up to the shear induced stress independently on
the shear amount $c$, i.e.:
\begin{equation}
  \label{eq:total-stress}
  \tens{\sigma}(\bu_{ij}) =  \tens{\sigma}(\bu_{ij}) + \tens{\sigma}_p
  = c \mu (\bv{e}_i\otimes\bv{e}_j + \bv{e}_j\otimes\bv{e}_i)
  +
  p \frac{(2\mu+\lambda)}{\mu} \tens{\beta}.
\end{equation}

This modified expression of the stress should be taken into account when performing pure
shear measurements in \emph{in-vivo} pressurized tissues, as the
measured forces will contain also the pressure-induced term
\begin{equation}
  \label{eq:delta-shear}
  F^p_{ij} = -|A| p \frac{(2\mu+\lambda)}{\mu} \sum_{k=1}^3
  \beta_k \btau_{k}\cdot \bv{e}_i \btau_{k}\cdot\bv{e}_j,
\end{equation}
which corresponds to the measured force due to vessels in the
direction $\bv{e}_j$, measured on the face with area $|A|$ and normal
$\bv{e}_i$.  
In practice, the effective shear modulus $\mu^e$ that
would be measured with a pure shear experiment with shear displacement
of scale $c$, is offset with respect to the shear modulus $\mu$ of the
elastic matrix by a factor that depends on the orientation of the
vessels, their volume fraction, the amount of internal vessel pressure,
and the applied shear displacement $c$:
\begin{equation}
  \label{eq:shear-modulus-correction}
  \mu^e_{ij} = \left(1- p  \frac{(2\mu+\lambda)}{c\mu^2} \sum_{k=1}^3
  \beta_k \btau_{k}\cdot \bv{e}_i \btau_{k}\cdot\bv{e}_j\right) \mu.
\end{equation}

Equation~\eqref{eq:shear-modulus-correction}
may explain some of the experimental observations in the
literature \cite{Chatelin2011}, where the in vitro shear modulus values
obtained by transient elastography were much lower than those found in vivo
(with a mean difference of 66\%).

\section{The discrete problem}\label{ssec:discrete}

Let us now assume to deal with a polygonal or polyhedral domain
$\Omega$, and let $\{\mathcal T_h\}_h$ be a family of conformal,
quasi-uniform quadrilateral or hexaedral meshes exactly covering
$\Omega$ where $h$ denotes the maximum element diameter.

We construct the finite element space
of globally continuous piecewise polynomials of
order $k$ in each coordinate directions defined by:
\begin{equation}
  \label{eq:definition-Vh}
  V_{h} 
  := \{ \vv \in (H^1_{\Gamma_D}(\Omega))^d, \text{ s.t. } \vv|_T \in
  \mathcal Q^k(T) \quad \forall T \in \mathcal T_h, \vv |_{\Gamma_D} =
  0\},
\end{equation}
denoting with $m$ its dimension and with $\{ \hat \vv_i\}_{i=1}^{m}$ a
basis for the space.

We consider the following discrete problem:
\begin{problem}[3D, Discrete] 
  \label{pb:discrete}
  Let be given a curve $\bgamma: [0,L] \to \mathbb R^3$, a function
  $a: [0,L] \to \mathbb R$ describing the radius, and a pressure
  $p_{\bgamma}: [0,L] \to \mathbb R$. Moreover, let $\bv{F}$ be one of
  the singular or hyper-singular source terms defined in the previous sections.  
  Find the
  displacement $\bu_h \in V_{h}$ such that
  \begin{equation}
    \label{eq:discrete}
          (2\mu \tens{e}(\bu_h), \tens{e}(\vv_h))_{\Omega} + (\lambda \nabla\cdot
          \bu_h, \ldiv\vv_h)_\Omega= <\bv{F},\vv_h>, \qquad \forall \vv_h
          \in V_{h}.
    \end{equation}
\end{problem}

Problem~\ref{pb:discrete} reduces to the solution of the following
linear system of equations
\begin{equation}
  \label{eq:linear-system}
\tens{K} \cdot \bv U = \bv b,
\end{equation}
where 
\begin{equation}
  \label{eq:definition-K}
  \begin{split}
    \bv K_{ij}  &:= 2\mu (\tens{e}(\hat \vv_j), \tens{e}(\hat \vv_i))_{\Omega} + \lambda
    (\ldiv\hat \vv_j, \ldiv\hat \vv_i)_\Omega\\
      \bv b_{i}  &:= <\bv{F}, \hat \vv_i>_\Omega,
  \end{split}
\end{equation}
and $\bv U = \{u^i\}_{i=1,\hdots,m}$ indicates the vector of coefficients of the finite element function 
$\bu_h$  such that
\begin{equation}
  \label{eq:definition-uh}
  \bv u_h(\bx) = \sum_{i=1}^m u^i \vv_i(\bx) \quad \bx \in \Omega.
\end{equation}

The right-hand side of equation \eqref{eq:discrete} contains a
singular forcing term, whose numerical computation may require a discrete
approximation of the Dirac delta distribution.
We use one of
the classical approximations widely employed in the context of the Immersed
Boundary Method~\cite{Peskin-2002-The-immersed-boundary-0}, i.e., 
\begin{equation}\label{eq:dirac-approximation}
\delta^\varepsilon(\bx)  := \frac{1}{\varepsilon^d} \prod_{i=1}^d \theta\left( \frac{\bx_i}{\varepsilon}\right)
\end{equation}
with
\begin{equation}
    \theta(y)  := 
    \begin{cases} 
      \left(\cos(\pi y)+1\right)/2 & \text{ if } -1 < y < 1 \\
      0 & \text{otherwise}
    \end{cases}
\end{equation}
where $\varepsilon$ is an arbitrary (small) parameter.
In the numerical experiments presented in the following Section,
we set $\varepsilon = 2h$, i.e., twice the diameter of the smallest triangulation element.

This approximation of the Dirac distribution guarantees that
$\int_{\Re^d} \delta^\varepsilon \d \bx = 1$ and that
$\delta^\varepsilon \in C^1(\Re^d)$, making it a good candidate for a
regularization of the Dirac distribution required for the
hyper-singular formulation.
For a in-depth discussion on alternative approximations of the Dirac distribution
and of their approximation
properties, we point the reader to the excellent
work of~\cite{Hosseini2014}.

\section{Numerical results}\label{sec:numerics}
This section is dedicated to the numerical validation of the mathematical models derived in Sections 
\ref{sec:2d} and \ref{sec:3d}. We will consider first the simple 2D axi-symmetric situation with known exact solution.
Next, we will investigate the effect of random distribution of vessels (in two and three dimensions) and
use the hyper-singular formulation to derive a statistical model for the effective tissue behavior.

All numerical examples provided in this section were obtained using an
open source code based on the \texttt{deal.II}
library~\cite{Bangerth2007,ArndtBangerthDavydov-2017-a,AlzettaArndtBangerth-2018-a}. The
code is freely available at the address
\url{https://gitlab.com/code_projects/immersed-elasticity}
(\cite{HeltaiCaiazzo-2018-a}), and it is inspired by the \texttt{deal.II} \emph{step-60}
tutorial~\cite{HeltaiAlzetta-2018-a}.
All simulations were performed using $\mathcal Q^1$ conforming finite
element spaces on quadrilaterals or hexaedral meshes.

\subsection{Reference solution in two-dimension}
We consider first a 2D axi-symmetric problem, comparing 
the results obtained via the proposed method employing
three different source terms to approximate the vessel network:
\begin{enumerate}
\item[{\sf{(S)}}] 
A singular forcing term, whose
  distributional definition is given by
  \begin{equation*}
    <\bv{F}, \vv> := \int_{\Gamma^a} \frac{(2\mu+\lambda)}{\mu} p \vv
    \d \Gamma
\end{equation*}
\item[ {\sf{(RS)}}] A regularized singular forcing term, given by
  \begin{equation*}
    \bv{F}^\varepsilon (\bx) := \int_{\Gamma^a} \frac{(2\mu+\lambda)}{\mu} p
    \delta^\varepsilon(\bx-\by) \d \Gamma_y
\end{equation*}
\item[ {\sf{(RHs)}}] A regularized hyper-singular forcing term, given by
  \begin{equation*}
    \bv{F}^\varepsilon (\bx) := -\frac{(2\mu+\lambda)}{\mu} \pi a^2 p 
    \grad \delta^\varepsilon(\bx) 
\end{equation*}
\end{enumerate}

\begin{figure}
   \centering
  \includegraphics[width=.32\textwidth]{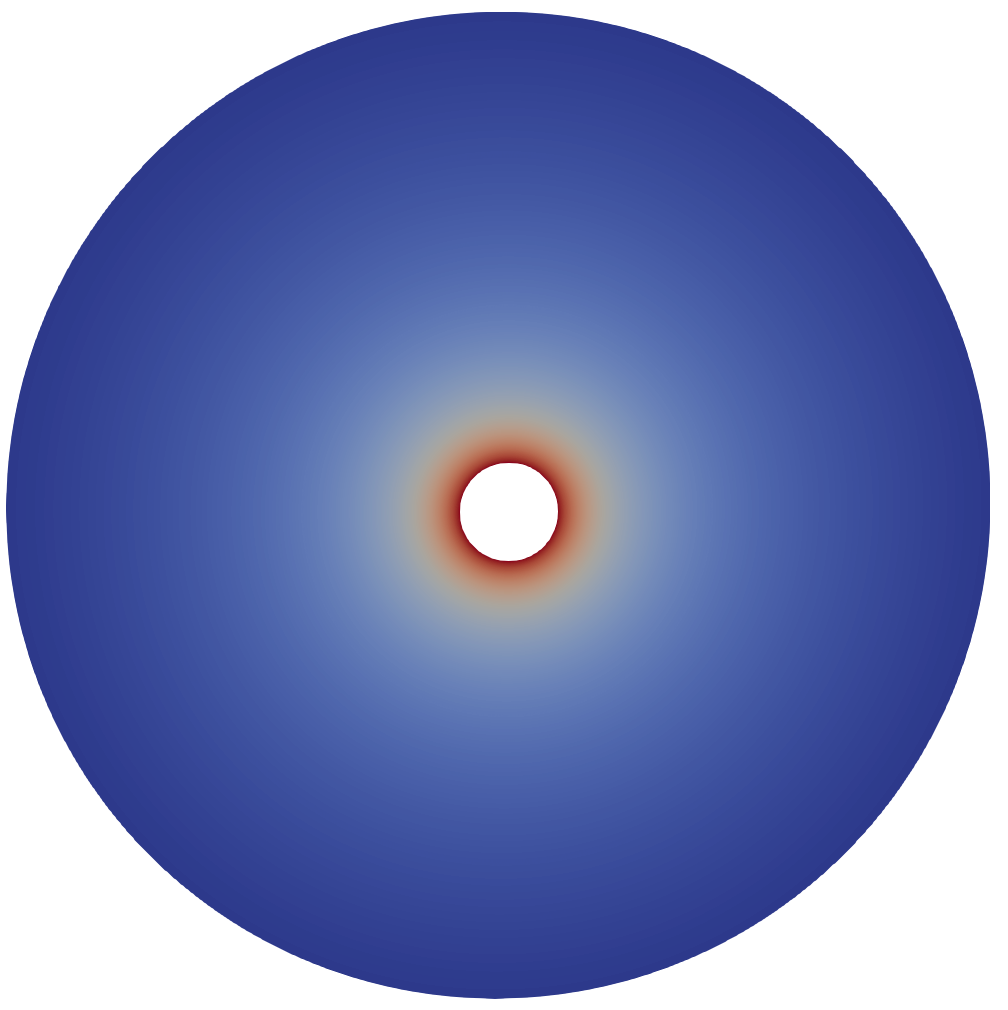}
  \hspace{1cm}
  \includegraphics[width=.32\textwidth]{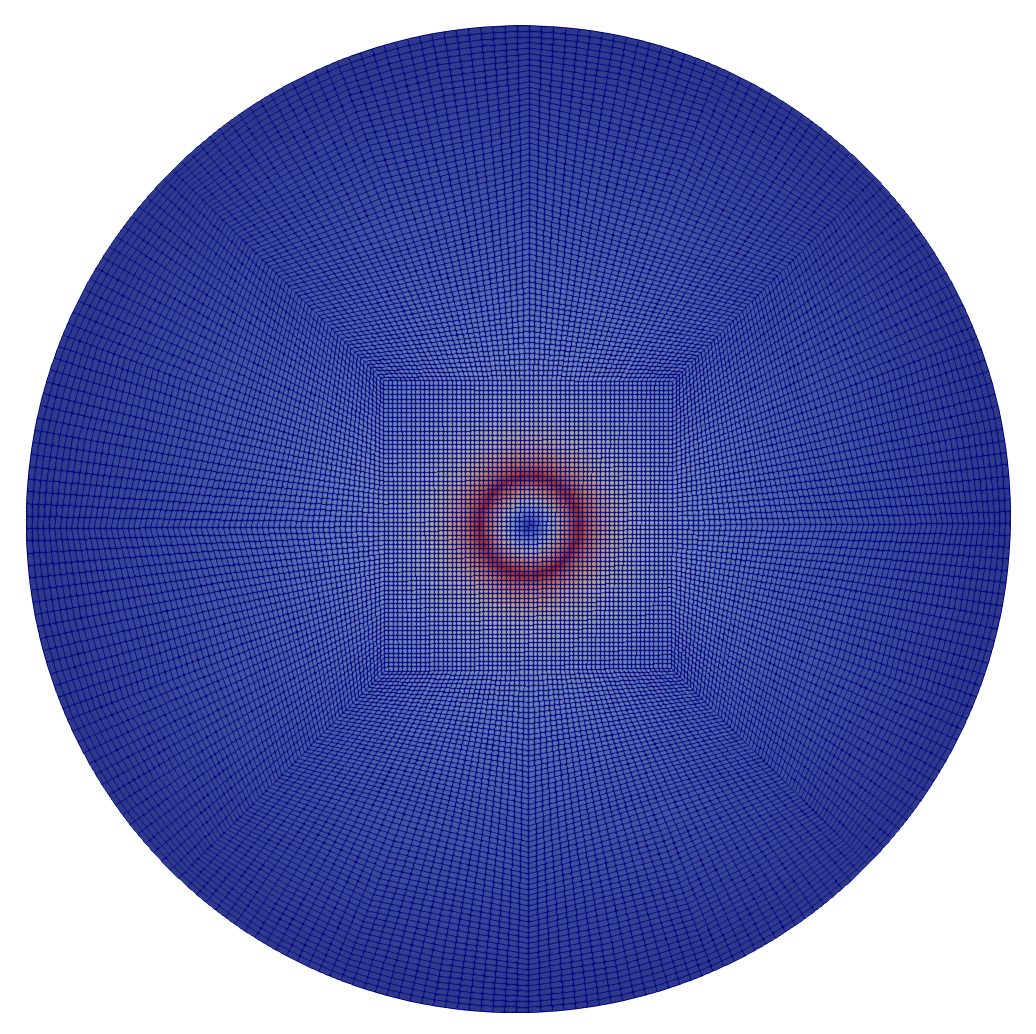}
  \caption{2D axi-symmetric problem. Comparison between the exact solution (left) and 
  the numerical solution obtained with the singular forcing term (case {\sf{(S)}}).}
  \label{fig:comparison-mesh}
\end{figure}

We consider $\lambda = \mu = p = 1$ Pa.
Figure~\ref{fig:comparison-mesh} 
shows the comparison,
on a circular domain of radius $R=1$,
with a vessel of radius $a=.1$, between the exact solution (left plot) and 
the solution obtained using the singular source  {\sf{(S)}} (right plot).

Notice that the grid does not need to be aligned to the surface
$\Gamma$, thanks to non-matching interpolation techniques~\cite{RoyHeltaiCostanzo-2015-a, HeltaiCostanzo-2012-a,
  BoffiGastaldiHeltaiPeskin-2008-a, Heltai-2008-a}.

 %
 Although the ratio $\frac{a}{R}$ is not too small (equal to 0.1), the
 solutions are remarkably close outside of the vessel, in agreement
 with the asymptotic expansion \eqref{eq:2d-sigma-jump-approx}, which predicts
 a residual error of the order of $\sim \mathcal{O}\big((a/R)^2\big)$.


 For a more quantitative assessment, we studied the error, with
 respect to the known exact solution, of the standard finite element
 approximation in the exact domain (applying a Neumann boundary
 condition on the resolved vessel boundary) and of the regularized
 hyper-singular approach {\sf{(RHs)}}.
We considered two situations, with vessel radius $a=0.1$ and $a=0.01$.

\begin{figure}
\centering
  \tikzsetnextfilename{H1_error_comparison_a.1}
  \begin{tikzpicture}
    \begin{loglogaxis}[
      width=.45\textwidth,
      height=.35\textwidth,
       xlabel={Degrees of freedom},
      ylabel={$\|u-u_h\|_{H^1}$ },
      grid=major,
      legend entries={{\small Exact domain}, {\small Immersed: Approximated Dirac}}
      ]
      \addplot table[x=1, y=4]{\DataTwoD};
      \addplot table[x=1, y=4]{\DataTwoDDirac};

      \logLogSlopeTriangleReversed{0.8}{0.3}{0.1}{.5}{black}{1.0};
      
    \end{loglogaxis}
  \end{tikzpicture}
  \hfill
  %
  \tikzsetnextfilename{L2_error_comparison_a.1}
  \begin{tikzpicture}
    \begin{loglogaxis}[
      width=.45\textwidth,
      height=.35\textwidth,
      xlabel={Degrees of freedom},
      ylabel={$\|u-u_h\|_{L^2}$ },
      grid=major,
      legend entries={{\small FEM Exact domain}, {\small Immersed: Approximated Dirac}}
      ]
      \addplot table[x=1, y=2]{\DataTwoD};
      \addplot table[x=1, y=2]{\DataTwoDDirac};

      \logLogSlopeTriangleReversed{0.8}{0.3}{0.1}{1.0}{black}{2.0};
      
    \end{loglogaxis}
  \end{tikzpicture}

  \caption{Comparisons of the errors with respect to the analytical solution obtained
  using a finite element method where the mesh resolved the vessel-tissue interface and using the
  immersed (approximated Dirac) approach, for
    vessel radius $a=.1$.
  Left. Error in $H^1$norm. Right. Error $L^2$-norm. The triangles refer to first order  (left plot) and second order (right plot) 
  convergence slopes. The plots show that the immersed approach has the same order of 
  accuracy, and a comparable numerical error, as the
  full finite element method, up to a mesh size comparable $O(a^2)$. }
  \label{fig:error_comparison_H1_a=.1}
\end{figure}
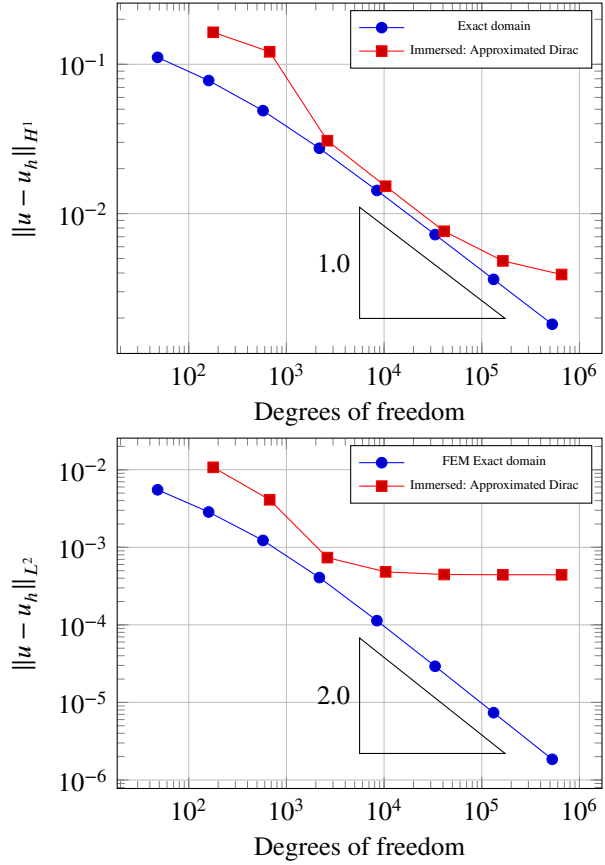

\begin{figure}
\centering
  \tikzsetnextfilename{H1_error_comparison_a.01}
  \begin{tikzpicture}
    \begin{loglogaxis}[
      width=.45\textwidth,
      height=.35\textwidth,
      xlabel={Degrees of freedom},
      ylabel={$\|u-u_h\|_{H^1}$ },
      grid=major,
      legend entries={{\small FEM Exact domain}, {\small Immersed: Approximated Dirac}}
      ]
      \addplot table[x=1, y=4]{\DataTwoDSmall};
      \addplot table[x=1, y=4]{\DataTwoDDiracSmall};

      \logLogSlopeTriangleReversed{0.8}{0.15}{0.1}{.5}{black}{1.0};
      
    \end{loglogaxis}
  \end{tikzpicture}
  \hfill
  %
  \tikzsetnextfilename{L2_error_comparison_a.01}
  \begin{tikzpicture}
    \begin{loglogaxis}[
      width=.45\textwidth,
      height=.35\textwidth,
      xlabel={Degrees of freedom},
      ylabel={$\|u-u_h\|_{L^2}$ },
      grid=major,
      legend entries={{\small FEM Exact domain}, {\small Immersed: Approximated Dirac}}
      ]
      \addplot table[x=1, y=2]{\DataTwoDSmall};
      \addplot table[x=1, y=2]{\DataTwoDDiracSmall};

      \logLogSlopeTriangleReversed{0.8}{0.15}{0.1}{1.0}{black}{2.0};
      
    \end{loglogaxis}
  \end{tikzpicture}
  \caption{Comparisons of the errors with respect to the analytical solution obtained
  using a finite element method where the mesh resolved the vessel-tissue interface and using the
  immersed (approximated Dirac) approach, for
    vessel radius $a=.01$.
  Left. Error in $H^1$norm. Right. Error $L^2$-norm. 
  The triangles refer to first order  (left plot) and second order (right plot) 
  convergence slopes. 
  The plots show that the immersed approach has the same order of 
  accuracy, and a comparable numerical error, as the
  full finite element method, up to a mesh size comparable $O(a^2)$. }
  \label{fig:error_comparison_H1_a=.01}
\end{figure}
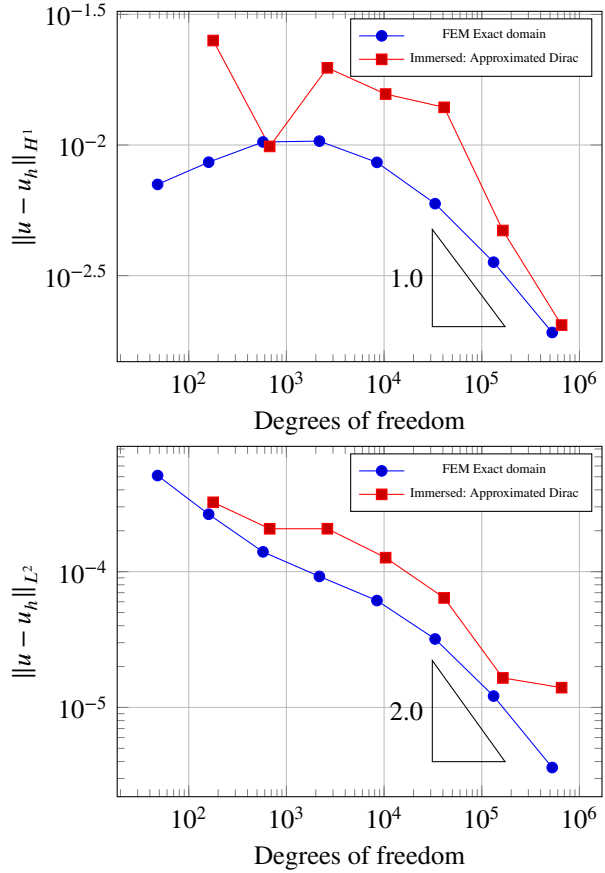

For the first case (radius $a=0.1$), the errors are reported in Tables \ref{tab:exact_hole_a=.1}
and~\ref{tab:approx_hole_a=.1}, and graphically in
Figure~\ref{fig:error_comparison_H1_a=.1}.
Tables~\ref{tab:exact_hole_a=.01} and~\ref{tab:approx_hole_a=.01}
report the errors for the exact domain case and the regularized
hyper-singular case when the vessel radius is $a=.01$, while, in this case, the graphical
comparison is reported in Figure~\ref{fig:error_comparison_H1_a=.01}.

\begin{table}
  \centering
  \input{tables/error_exact_hole_a.1.tex}
  \caption{Error on the standard Finite Element Approximation, with
    exact domain, and vessel radius $a= 0.1$. The columns reporting the number of degrees of freedom
    is particularly relevant for the comparisons, in terms of efficiency, with the immersed formulation.}
  \label{tab:exact_hole_a=.1}
\end{table}

\begin{table}
  \centering
  \input{tables/error_approx_dirac_a.1.tex}
  \caption{Error on the Finite Element Approximation, with
    approximated vessel using the regularized hyper-singular approach
    \sf{(RHs)}, and vessel radius $a=0 .1$. The columns reporting the number of degrees of freedom
    is particularly relevant for the comparisons, in terms of efficiency, with the full finite element formulation.}
  \label{tab:approx_hole_a=.1}
\end{table}


%
%
%

The comparison clearly shows that for large vessels the exact domain
approximation remains the method of choice in terms of accuracy per
degree of freedom. 
In fact, in the case of large vessel radius, the generation of the discrete domain
does not increase the complexity of the overall simulation and does not increase
substantially the required resolution of the mesh.
As it can be expected, the simulation based on the exact domain has a clear advantage
with respect to the regularized approach only when the mesh size decreases below $a^2$, i.e., the square of vessel radius.
At this stage, the immersed method
exhibits a
plateau on the $L^2$ error, coherently with the asymptotic analysis
presented in the previous sections.

\begin{table}
  \centering
  \input{tables/error_exact_hole_a.01.tex}
  \caption{Error on the standard Finite Element Approximation, with
    exact domain, and vessel radius $a= .01$.}
  \label{tab:exact_hole_a=.01}
\end{table}

\begin{table}
  \centering
  \input{tables/error_approx_dirac_a.01.tex}
  \caption{Error on the Finite Element Approximation, with
    approximated vessel using the regularized hyper-singular approach
    \sf{(RHs)}, and vessel radius $a= .01$.}
  \label{tab:approx_hole_a=.01}
\end{table}

When the vessel size decreases, however, the domain generation
becomes more and more computationally expensive, the required
mesh characteristic sizes decreases considerably, and mesh
quality might deteriorate, challenging also the exact domain approach, as
shown in Figures~\ref{fig:error_comparison_H1_a=.01}. 
Hence, for small vessel radii, it has to be expected that full numerical simulations with mesh size 
below the critical $a^2$-regime are no longer feasible.

These last plots show that the correct order of convergence is reached
by the Finite Element Approximation on the fully resolved domain only using a large number of degrees of
freedom. On the other hand, for smaller vessel size, the
regularized hyper singular approach achieves a comparable accuracy
with respect to the exact domain approach, at a fraction of the
computational cost in terms of mesh generation, which is totally
independent on to the vessel geometry and location.

\begin{remark}
Although the accuracy of the fully resolved approach might still be better, the order
of magnitude of the errors are comparable for the two methods, 
making
the regularized approach a competitive alternative when the
ratio $a/R$ is small, with the
additional advantage that it does not require the meshing of the tissue-vessel interface,
which might introduce additional complexity in the fully resolved case.
\end{remark}


%
%
%

\subsection{Prototypical vessel junction}
The simple axis symmetric two-dimensional cases presented in the previous example provides
evidence that, when an exact solution is available, the proposed method achieves 
the correct order of convergence.

\begin{figure}[!h]
  \centering
  \includegraphics{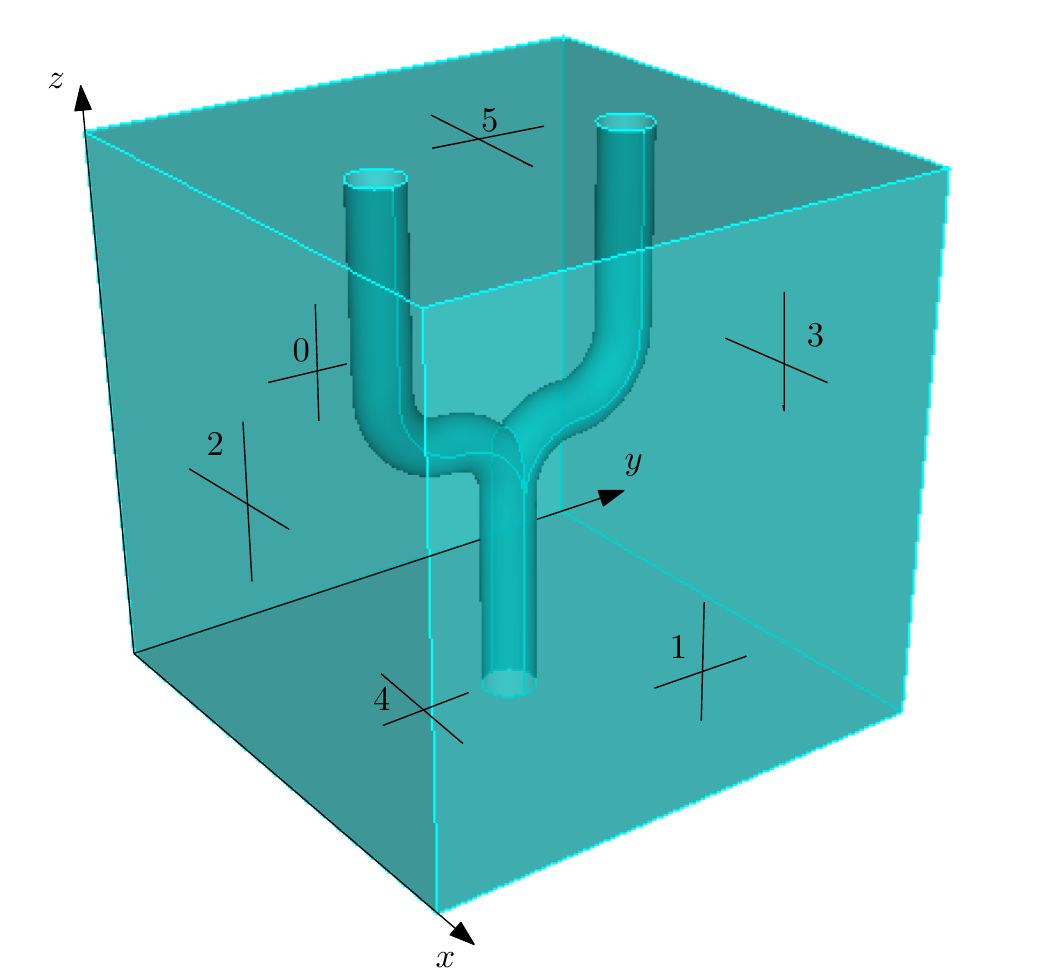}
  \caption{Prototypical building block of blood vessels: a
    $Y$-junction, where a single vessel splits into two. The
    $Y$-junction is surrounded by an elastic cube of dimension one,
    while the diameter of the vessel remains constant, and it is equal
    to $0.1$. The minimum radius of curvature of the vessel centerline
    is equal to $0.125$. Numbers on the faces are used in
    Table~\ref{tab:comparison-y-junction} to refer to the average
    displacements. }
  \label{fig:branching-test}
\end{figure}

To further assess the prediction capabilities of our method, we consider 
a case in three dimensions that represents the prototypical building
block of blood vasculature: a $Y$-junction branching (see
Figure~\ref{fig:branching-test}). 
We consider a $Y$-junction pressurised vessel of diameter $0.1$,
passing through a cubic sample of dimensions $[0,1]^3$ , clamped at
the top and at the bottom, and stress free on the lateral
surfaces. The minimum radius of curvature of the vessel centerline is
equal to $0.125$.

The setting is designed to understand
(i) the effect of the  curvature
of the vessel, verifying that the immersed approximation
is able to deliver a good approximation at the effective tissue scale, and 
(ii) verify that
the model works when vessels intersect and overlap.

In  this case, we compare the results obtained using with a finite element simulation based on a 
fine mesh on the exact domain (i.e., resolving the vessel-tissue interafce), with the results obtained with
the regularized hyper-singular forcing term \eqref{eq:3d-hyper}:
  \begin{equation*}
    \bv{F}^\varepsilon (\bx) := -\int_\Gamma \frac{(2\mu+\lambda)}{\mu} \pi a^2 p 
    \grad \delta^\varepsilon(\bx).
\end{equation*}
In the case of the immersed method, the discrete domain does not need to resolve the vessel, but the mesh is 
adaptively refined near the junction centerline.
We set $\lambda = \mu = p = 1$ Pa.

Table~\ref{tab:comparison-y-junction} provides a comparison of the
average displacements on the lateral and front faces:
\begin{table}[!h]
  \centering
  \begin{tabular}{lc|c}
    & Exact domain & Hypersingular\\
   $\bu_0 $ & 
 $\begin{pmatrix}
  -3.87e-03\\
  1.48e-06\\
  -3.41e-04\\
\end{pmatrix}$  & 
 $\begin{pmatrix}
  -4.81e-03\\
  -4.66e-13\\
  -1.75e-04\\
\end{pmatrix}$  \\
$\bu_1 $ & 
 $\begin{pmatrix}
  3.87e-03\\
  1.48e-06\\
  -3.41e-04\\
\end{pmatrix}$  & 
 $\begin{pmatrix}
  4.81e-03\\
  -4.25e-13\\
  -1.75e-04\\
\end{pmatrix}$  \\
$\bu_2 $ & 
 $\begin{pmatrix}
  -8.86e-11\\
  -3.18e-03\\
  -7.70e-04\\
\end{pmatrix}$  & 
 $\begin{pmatrix}
  3.73e-13\\
  -2.88e-03\\
  -2.07e-04\\
\end{pmatrix}$  \\
$\bu_3 $ & 
 $\begin{pmatrix}
  -9.65e-11\\
  3.18e-03\\
  -7.69e-04\\
\end{pmatrix}$  & 
 $\begin{pmatrix}
  -2.93e-13\\
  2.88e-03\\
  -2.07e-04\\
\end{pmatrix}$
  \end{tabular}
  \caption{Comparison between exact domain and hyper-singular lateral
    displacements of the $Y$-junction problem. The numbering of the
    displacements follows the convention used in
    Figure~\ref{fig:branching-test}. In all cases, the leading order component of the average displacements obtained 
    with the full finite element simulation and with the hypersingular formulation are very similar.}
  \label{tab:comparison-y-junction}
\end{table}

The table shows that the two simulations agree both qualitatively and
quantitatively, providing a good estimation of the behaviour of the
tissue surrounding the junction, with average displacements in the
same order of magnitude. Given that the radius of the vessel is
$R=0.1$, the expansion used to derive the hypersingular model, expects
an error in the average displacements in the order $10^{-2}$, i.e.,
$R^2$. The simulations, however, show a quantitative agreement with
maximum errors of one order of magnitude less than $R^2$.

\begin{figure}
   \centering
  \includegraphics[width=.45\textwidth]{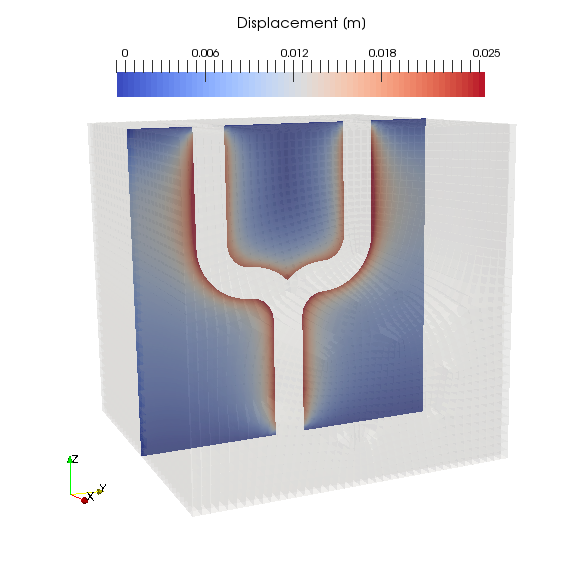}
  \hfill
   \includegraphics[width=.45\textwidth]{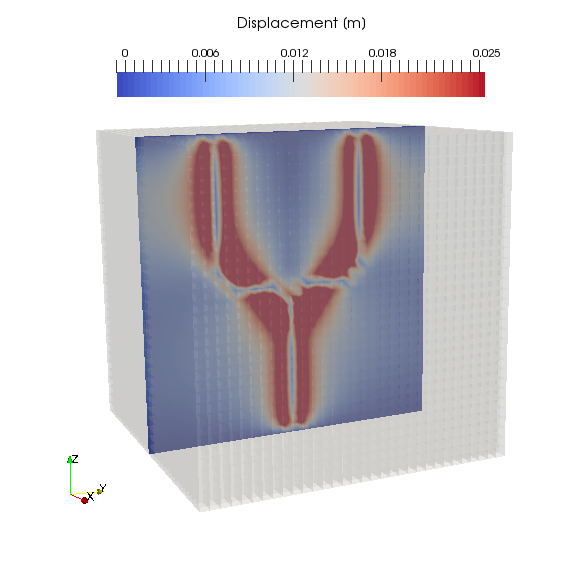}\\
  \includegraphics[width=.45\textwidth]{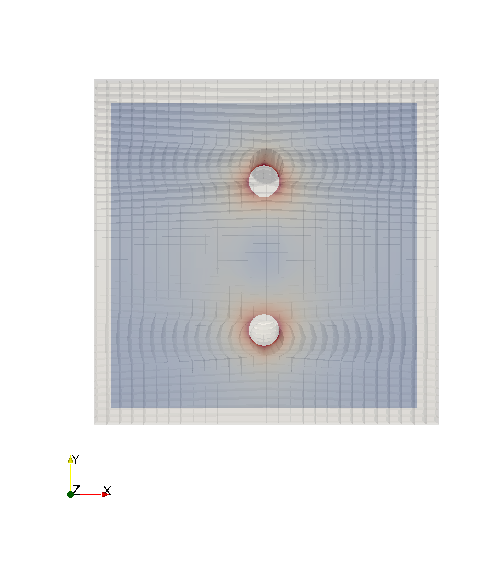}
  \hfill
   \includegraphics[width=.45\textwidth]{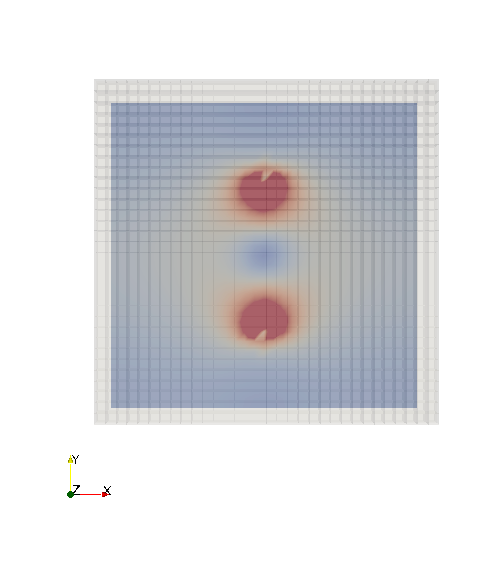}
   \caption{Local displacement plots in the exact domain case (left
     figures), and in the hyper-singular model (right figures). Away
     from areas of high curvature (bottom figures), the local
     agreement is good. Close to areas of high curvature, the
     hyper-singular model tends to flatten the response of the elastic
     material.}
  \label{fig:junction}
\end{figure}

Inspection of the displacement field close to the vessels is provided
in Figure~\ref{fig:junction}. The displacement plots show that there
is good agreement between the solutions in planes that are
perpendicular to the vessel direction, far from regions of high
curvature of the vessels (bottom figures).

However, the hyper-singular model tends to flatten the response of the
material around areas with very high curvature (top right figure),
providing an effective response similar to the one coming
from a \emph{straightened} version of the vessels. While the local
response is clearly different in the two cases,
Table~\ref{tab:comparison-y-junction} shows that the average response
on the faces of the sample is within the expected range, and the hypersingular
formulation yields, on average, comparable displacement in the leading order of magnitude.

A correction of the model that takes explicitly into account the
curvature of the vessels is currently under investigation.

\subsection{Pressurized tissue in two dimensions}

\begin{figure}[!h]
   \centering
 \includegraphics[width=.42\textwidth,trim=2cm 1cm 7cm 0cm,clip=true]{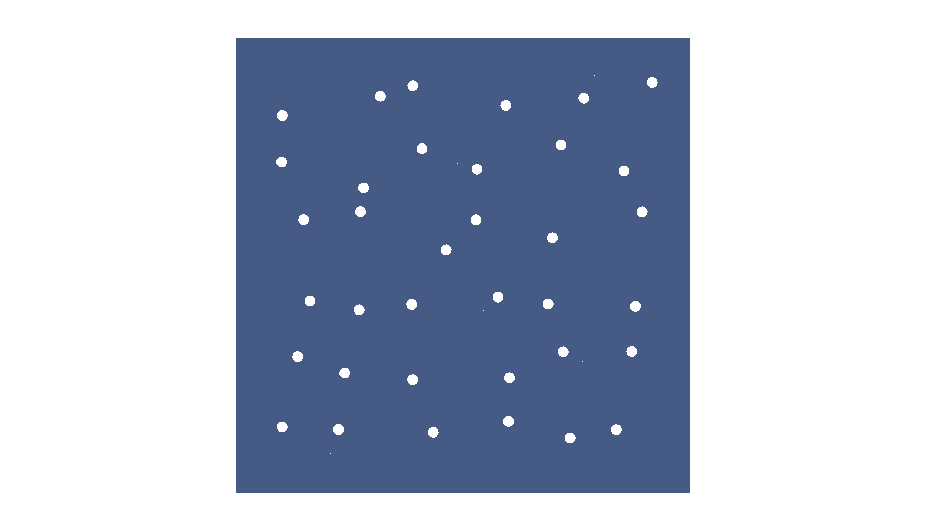}
  \hspace{.5cm}
   \includegraphics[width=.25\textwidth,trim=6.5cm 3cm 7.5cm 2cm,clip=true]{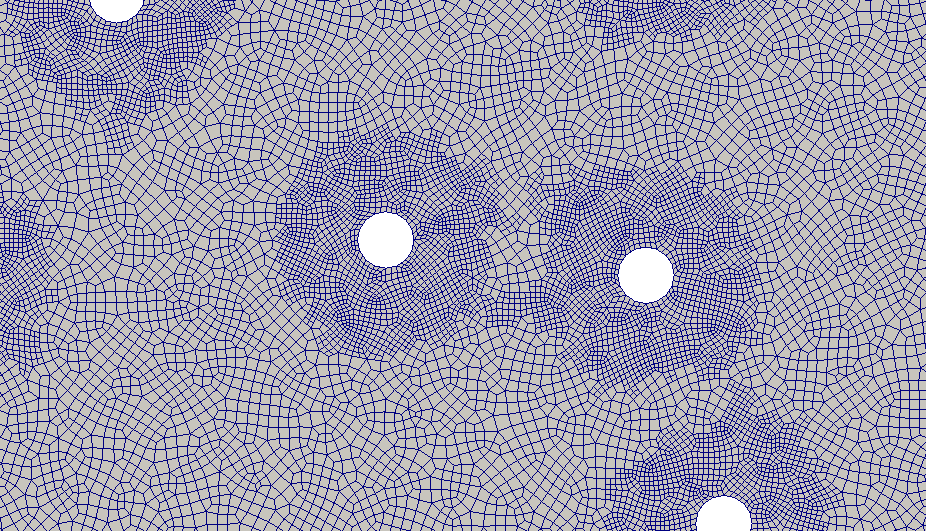}
  \hspace{.1cm}
  \includegraphics[width=.25\textwidth,trim=6.5cm 3cm 7.5cm 2cm,clip=true]{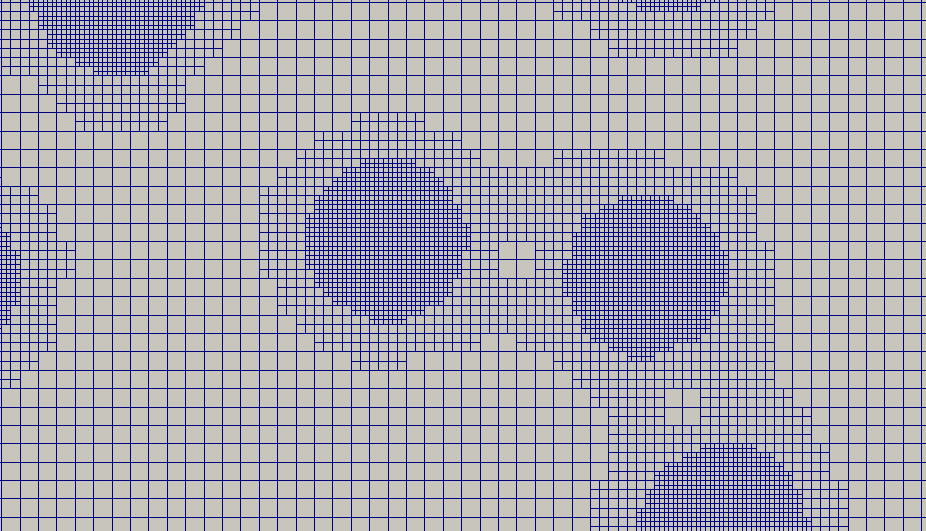}
  \caption{Left: Two-dimensional domain (unit square) with 36 vessels of radius $r=0.12$.
  Right: Mesh around a vessel in the full-scale (left) and immersed (right) simulation.}
  \label{fig:tissue2d-domain}
\end{figure}

As next, we consider a tissue sample (2D) $[0,1] \times [0,1]$ with a set of 36 vessels of
radius $r=0.012$ placed at random locations (see Figure \ref{fig:tissue2d-domain}, left).
In each vessel, we prescribe an unitary pressure $p=1$, comparing the results
for the tissue displacement when performing a \textit{full-scale} simulation
(e.g., resolving the vessel interface and applying a Neumann boundary condition)
and when using the regularized immersed method. 

We consider the case where the immersed method has a computational
complexity (in terms of number of degrees of freedom) comparable with
the full-scale simulation.  Notice, however, that the immersed method
does not require the explicit resolution of the interface, hence
allowing for a much easier generation of the computational mesh, which
has been constructed starting from an uniform mesh on the unit square
via an automatic refinement strategy based on the Kelly-error
estimator (see, for example,~\cite{Kelly1983}).

Figure \ref{fig:tissue2d-domain} (right) shows the size of the discretization around the vessels
in both cases.

\begin{figure}[!h]
   \centering
  \includegraphics[width=.35\textwidth,trim=5.5cm 0cm 7.5cm 0cm,clip=true]{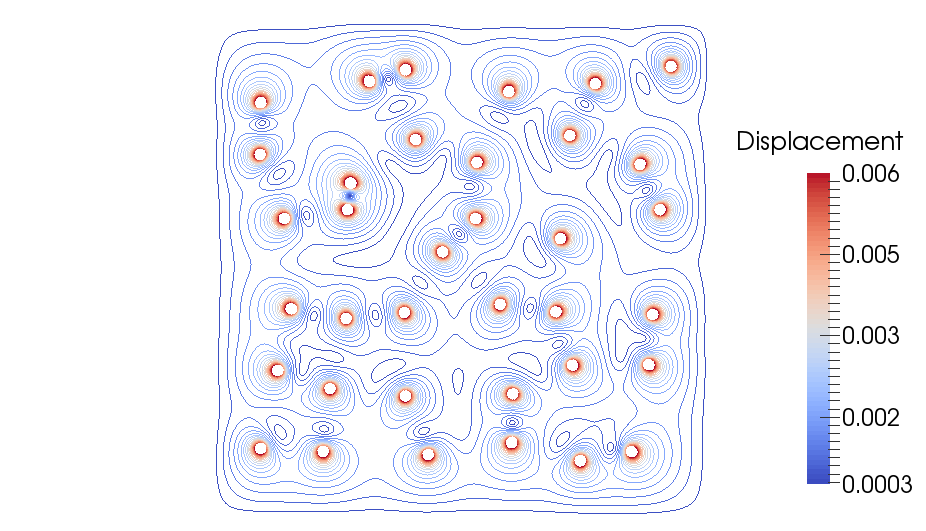}
  \hspace{.5cm}
  \includegraphics[width=.35\textwidth,trim=5.5cm 0cm 7.5cm 0cm,clip=true]{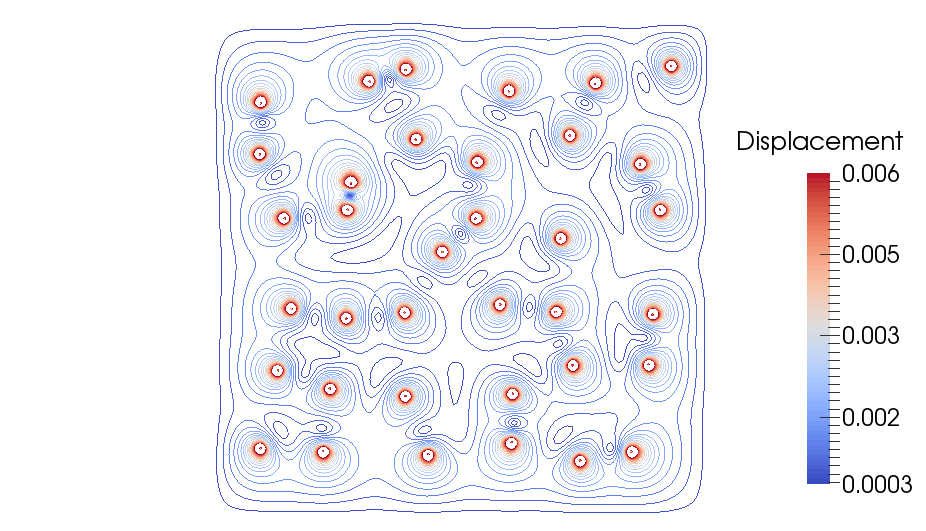}
  \hspace{1cm}
  \includegraphics[width=.17\textwidth,trim=25cm 0cm 0cm 0cm,clip=true]{figures/tissue2d_immersed_n36_contour.png}
  \caption{Two-dimensional pressurized tissue: isolines of the displacement field. Left: full-scale simulation.
  Right: Immersed method.}
  \label{fig:tissue2d-results}
\end{figure}

\begin{figure}[!h]
   \centering
  \includegraphics[width=.35\textwidth,trim=5cm 2cm 8.2cm 2cm,clip=true]{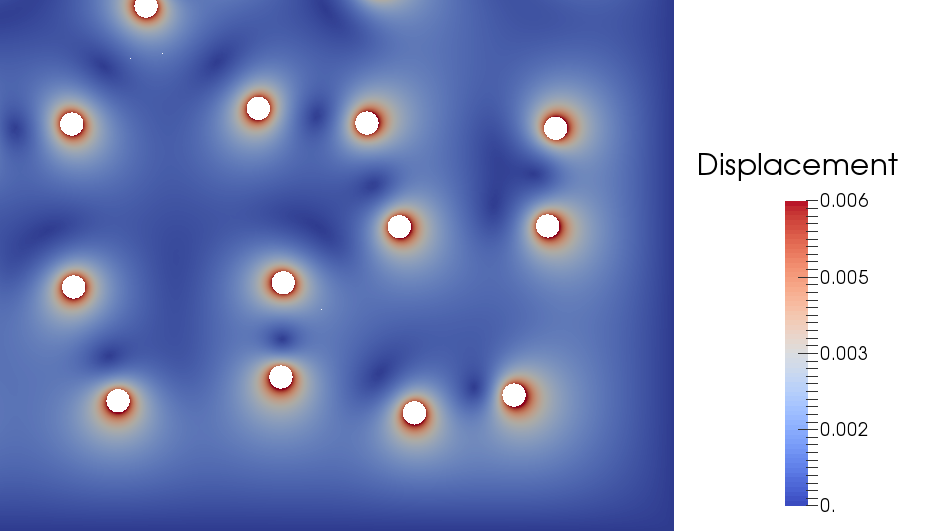}
  \hspace{.5cm}
  \includegraphics[width=.35\textwidth,trim=5cm 2cm 8.2cm 2cm,clip=true]{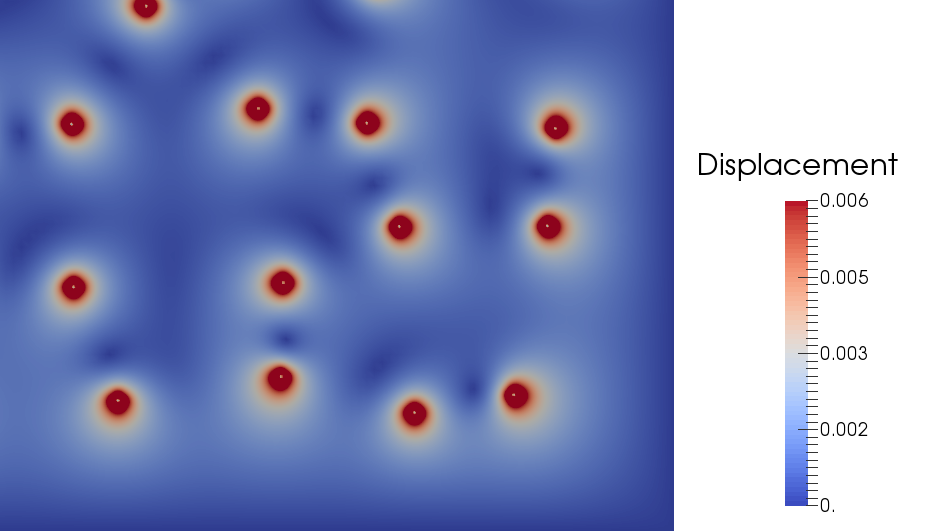}
  \hspace{1cm}
  \includegraphics[width=.17\textwidth,trim=24.5cm 2cm 0cm 0cm,clip=true]{figures/tissue2d_immersed_n36_zoom.png}
  \caption{Two-dimensional pressurized tissue: zoom of the numerical solution 
  (magnitude of the displacement field) close to the bottom-right corner, in order to better compare the values of displacement
  near the pressurized vessels. Left: full-scale simulation.
  Right: Immersed method.}
  \label{fig:tissue2d-results-zoom}
\end{figure}

A comparison of the two approaches is provided in Figures \ref{fig:tissue2d-results}-\ref{fig:tissue2d-results-zoom}, 
demonstrating the qualitative agreement between the displacement fields, especially
close to the vessels.

\subsubsection{Homogenized -- Two-dimensional case}
\label{sec:two-dimensional-case}
In this section, we probe the hypothesis derived in Section
\ref{ssec:derivation} by performing a set of statistical simulations
in two dimensions, in which we consider several realizations of a random
collection of vessels in a box domain $\Omega = [0,1]^2$, we impose
zero Dirichlet conditions, and we measure the effects of the
deformations on the boundary, by averaging the forces exerted by the
expanding solid on the faces.

In this set of tests, we fix $\mu=1$, $p=1$, and we vary $\lambda$,
the number of vessels, the radius, and the refinement level of the
grid, to understand the robustness of the method with respect to grid
size, vessel density, and material properties. In the presented simulation,
the grid size is given by $h=2^{-\text{ref}}$, where $\text{ref}$ will be denoted as
the refinement level. 

\begin{figure}
  \centering
  \tikzsetnextfilename{random_distribution}
  \begin{tikzpicture}
    \begin{axis}[%
      axis equal image,
      xmin = 0,
      xmax = 1,
      ymin = 0,
      ymax =1,
      yticklabels={},
      xticklabels={},
      ticks=none,
      ]
      
      \pgfplotstableread{./data/statistical_2d_000_hyper_singularities.gpl}\mytable;

      \addplot [
      only marks,
      mark size=1pt,
      mark options={black},
      ]
      table
      {\mytable};
    \end{axis}
  \end{tikzpicture}
  \hspace{1cm}
  \includegraphics[width=5.7cm]{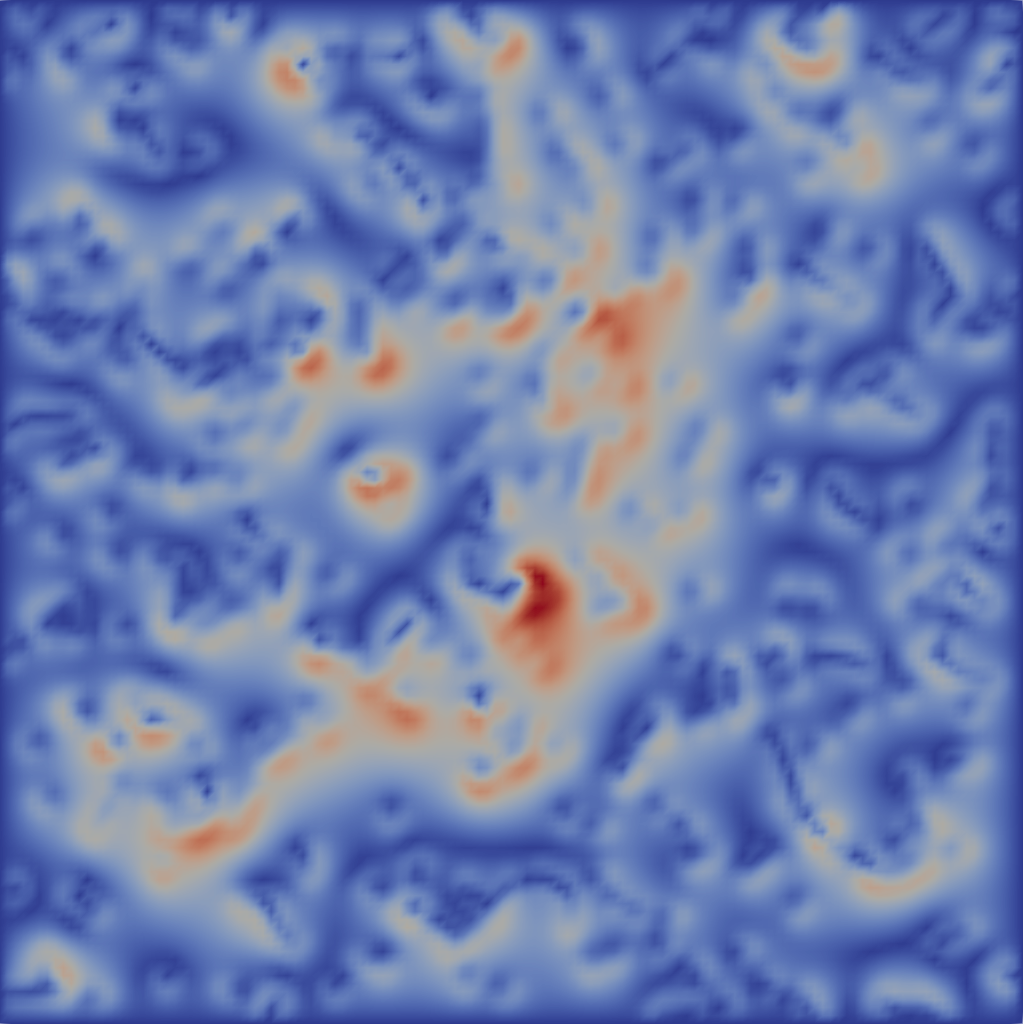}
    \caption{Left: Example of random realization with 500 vessels of radius
    $5.65\times 10^{-4}$ ($\beta=5\%$).  Displacement field, obtained with $\lambda = 1$
and mesh refinement level 7 ($h=2^{-7}  = 7.8125\times 10^{-3}$,
number of degrees of freedom equal to $33,282$). The maximum
    displacement is on the order of $10^{-3}$.}
    \label{fig:statistical-2d-500}
\end{figure}

The vessels are randomly distributed on $\Omega$, assuring that they
do not intersect the boundary.
An example of resulting distribution, considering $500$ randomly distributed
vessels of radius $5.65\times 10^{-3}$ (total volume fraction $\beta \approx 0.05$)
is provided in Figure~\ref{fig:statistical-2d-500} (left), while
Figure ~\ref{fig:statistical-2d-500} (right) shows the results 
obtained with $\lambda=1$.

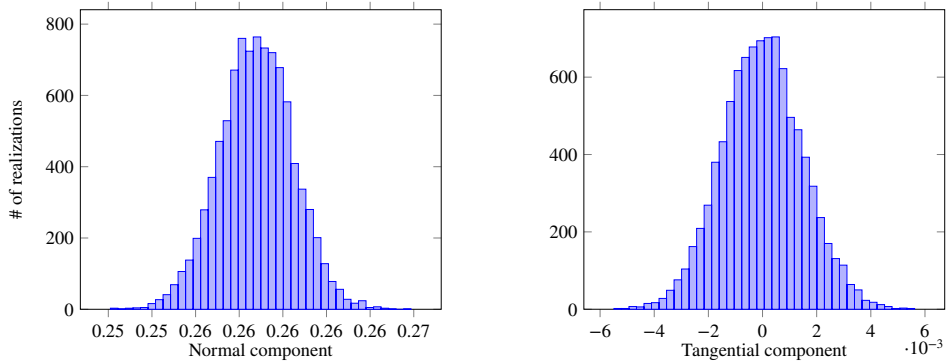
\begin{figure}
  \centering
  \tikzsetnextfilename{statistical_2d_normal}
  \begin{tikzpicture}[scale=.7]
    \begin{axis}[
      ybar,
      ymin=0,
      ylabel={\# of realizations},
      xlabel={Normal component}
      ]
      \addplot +[
      hist={
        bins=40,
      }   
      ] table [y index=0]
      {./data/statistical_2d,ref=9,lambda=1,n=30000,r=0.001.gpl};
    \end{axis}
  \end{tikzpicture}
  \hspace{1cm}
    \tikzsetnextfilename{statistical_2d_tangential}
  \begin{tikzpicture}[scale=.7]
    \begin{axis}[
      ybar,
      ymin=0,
      xlabel={Tangential component}
      ]
      \addplot +[
      hist={
        bins=40,
      }   
      ] table [y index=1]
      {data/statistical_2d,ref=9,lambda=1,n=30000,r=0.001.gpl};
    \end{axis}
  \end{tikzpicture}
  \caption{Statistical distribution of the normal (left) and
    tangential (right) pressure force, generated with
    $\beta \simeq 5\%$, $\lambda = 1$, $p=1$, and refinement
    refinement level 9 on 10,000 realizations.}
  \label{fig:statistical_distribution_2d}
\end{figure}

Figure~\ref{fig:statistical_distribution_2d} shows the statistical
distribution of the normal and tangential pressure force for 10,000
realizations of randomly distributed vessels with fixed volume
fraction $\beta \simeq 5\%$, Lam\`e parameter $\lambda = 1$, and
constant pressure $p=1$, on a grid with refinement refinement level 9
($h=2^{-9}=1.953125 \times 10^{-3}$).


\begin{figure}
 \centering
 \tikzsetnextfilename{statistical_distribution_normal_forces_ref8}
 \begin{tikzpicture}[scale=.7]
   \begin{axis}[%
     legend pos= north west,
     xtick={
       0.05, 0.1, 0.15
     },
     xlabel={Volume fraction of vessels $\beta$},
     ylabel={Normal component of $F$},
     ]
     \pgfplotstableread{./data/statistical_2d_error_bars_ref8.gpl}\mytable;
     \addplot+[
     only marks,
     mark size=1pt,
     error bars/.cd,
     y dir=both,
     y explicit
     ] table[x index = 0, y index = 1, y error index = 2]
     {\mytable};

    \addplot[only marks, mark size=1pt, mark
    options={red}] table[x index = 0, y index = 3]
    {\mytable};
     
    \legend{{Refinement level 8}, {Estimated}};
   \end{axis}
 \end{tikzpicture}
 \hspace{1cm}
 \tikzsetnextfilename{statistical_distribution_normal_forces_ref9}
 \begin{tikzpicture}[scale=.7]
   \begin{axis}[%
     legend pos= north west,
     xtick={
       0.05, 0.1, 0.15
     },
     xlabel={Volume fraction of vessels $\beta$},
     ]
     \pgfplotstableread{./data/statistical_2d_error_bars_ref9.gpl}\mytable;
     \addplot+[
     only marks,
     mark size=1pt,
     error bars/.cd,
     y dir=both,
     y explicit
     ] table[x index = 0, y index = 1, y error index = 2]
     {\mytable};

    \addplot[only marks, mark size=1pt, mark
    options={red}] table[x index = 0, y index = 3]
    {\mytable};

    \legend{{Refinement level 9}, {Estimated}};
     
   \end{axis}
 \end{tikzpicture}
  \caption{Mean and variance (vertical bars) value of the normal
    component of the force on the pressurized tissue computed from the
    statistical simulations, with refinement level 8 (left) and 9
    (right), compared with the value estimated in
    \eqref{eq:estimated-total-pressure-force}, for $\lambda = 1$
    (bottom three experimental measures) and $\lambda = 10$ (top three
    experimental measures).}
        \label{fig:force_vs_beta}
\end{figure}
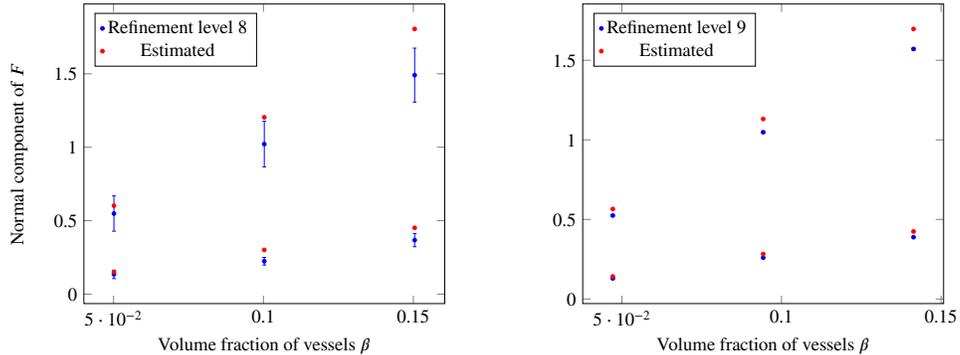

Figure \ref{fig:force_vs_beta} shows the mean normal force (with
$\lambda=1$ and $\lambda=10$) as a function of the volume fraction
$\beta$. The error bars show the value of the standard deviation,
which are only visible in the plot for the refinement 8 case, as in
the case with refinement 9 they are much smaller in scale. 

\subsubsection{Homogenized -- Three dimensional case (aligned vessels)}

In the three-dimensional setting, the homogenized hypothesis derived
in Section~\ref{ssec:derivation-3d} are more difficult to probe because
vessels are anisotropic in nature, and realizations with totally
random distributions of vessels cannot be reasonably considered as
uniformly distributed in the three directions \emph{at the same time},
making the homogeneized model only valid as an average approximated
model.

For simpler settings, where the vessels are uniformly distributed and
alligned along a preferred direction, the estimated forces of the
two-dimensional case presented in
Equation~\eqref{eq:estimated-total-pressure-force} are still a
reasonable estimate of the force generated by the pressurized vessels
measured on boundary walls whose normals are orthogonal to the vessels
center-line, but cannot be used to estimate the forces on walls whose
normal is parallel to the vessels. The two-dimensional setting
corresponds to an infinite material along the z-direction, where
deformation and stress are negligible along the z-direction.

\begin{figure}
  \centering
  \includegraphics[width=.4\textwidth]{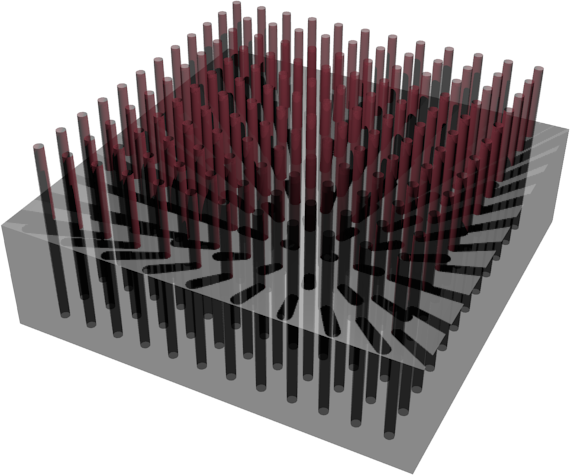}
  \caption{Schematic cut-view of a three-dimensional realization of a
    uniform distribution of vessels aligned along the z-axis.}
  \label{fig:3d-scheme}
\end{figure}

In the finite-domain case, boundary effects become more and more
important, and introduce a distortion in the homogenized estimate
given by Equation~\eqref{eq:estimated-total-pressure-force}. In this
section we consider several realizations of random collection of
streight vessels, aligned in the z-direction, included a box domain
$\Omega = [0, 1]^3$, as depicted in Figure~\ref{fig:3d-scheme}.

\begin{figure}
  \centering
  \tikzsetnextfilename{statistical_3d_lateral_normal}
  \begin{tikzpicture}[scale=.7]
    \begin{axis}[
      ybar,
      ymin=0,
      ylabel={\# of realizations},
      xlabel={Lateral normal component},
      xticklabel style={/pgf/number format/fixed,
        /pgf/number format/precision=4},
      ]
      \addplot +[
      hist={
        bins=40,
      }   
      ] table [y index=0]
      {./data/statistical_3d,ref=6,lambda=1,n=500,r=0.00565_sides.gpl};
    \end{axis}
  \end{tikzpicture}
  \hspace{1cm}    
    \tikzsetnextfilename{statistical_3d_lateral_tangential}
  \begin{tikzpicture}[scale=.7]
    \begin{axis}[
      ybar,
      ymin=0,
      xlabel={Lateral non-zero tangential component},
      xticklabel style={/pgf/number format/fixed,
        /pgf/number format/precision=4},
      ]
      \addplot +[
      hist={
        bins=40,
      }   
      ] table [y index=1]
      {./data/statistical_3d,ref=6,lambda=1,n=500,r=0.00565_sides.gpl};
    \end{axis}
  \end{tikzpicture}
  \caption{Statistical distribution of the lateral normal (left) and
    of the lateral non-zero tangential (right) pressure force,
    generated with $\beta \simeq 5\%$, $\lambda = 1$, $p=1$, and
    refinement refinement level 6 on 6,000 realizations.}
  \label{fig:statistical_distribution_3d_lateral}
\end{figure}
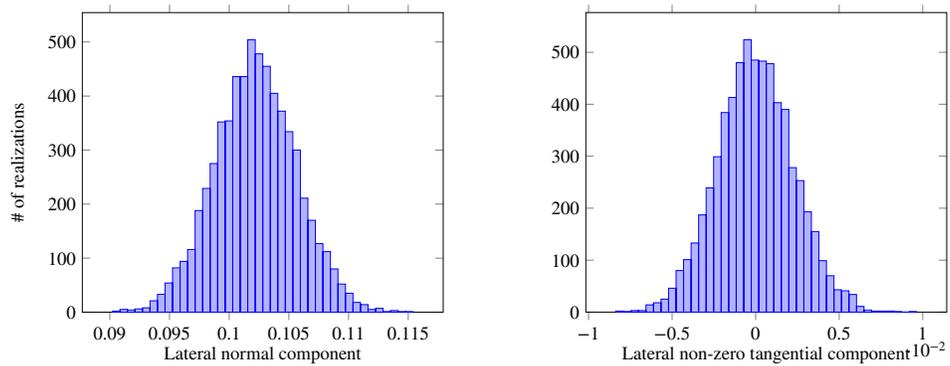

We impose zero Dirichlet conditions, and we measure the effect of the
deformations on the boundary due to the pressurized vessels, by
averaging the forces exerted by the expanding solid on the faces. By
symmetry, the z-component of the total force on each lateral face is
zero. In Figure~\ref{fig:statistical_distribution_3d_lateral}, we show
the statistical distribution of the non-zero components of the total
force on the lateral faces of the cube. 

\begin{figure}
  \centering
  \tikzsetnextfilename{statistical_3d_z_normal}
  \begin{tikzpicture}[scale=.7]
    \begin{axis}[
      ybar,
      ymin=0,
      ylabel={\# of realizations},
      xlabel={Top normal component},
      xticklabel style={/pgf/number format/fixed,
        /pgf/number format/precision=4},
      ]
      \addplot +[
      hist={
        bins=40,
      }   
      ] table [y index=0]
      {./data/statistical_3d,ref=6,lambda=1,n=500,r=0.00565_top.gpl};
    \end{axis}
  \end{tikzpicture}
  \hspace{1cm}
    \tikzsetnextfilename{statistical_3d_z_tangential}
  \begin{tikzpicture}[scale=.7]
    \begin{axis}[
      ybar,
      ymin=0,
      xlabel={Top $x$-component},
      xticklabel style={/pgf/number format/fixed,
        /pgf/number format/precision=4},
      ]
      \addplot +[
      hist={
        bins=40,
      }   
      ] table [y index=1]
      {./data/statistical_3d,ref=6,lambda=1,n=500,r=0.00565_top.gpl};
    \end{axis}
  \end{tikzpicture}
  \caption{Statistical distribution of the top and bottom normal
    (left) and along the $x$-direction (right) pressure force,
    generated with $\beta \simeq 5\%$, $\lambda = 1$, $p=1$, and
    refinement refinement level 6 on 1,500 realizations.}
  \label{fig:statistical_distribution_3d_z}
\end{figure}
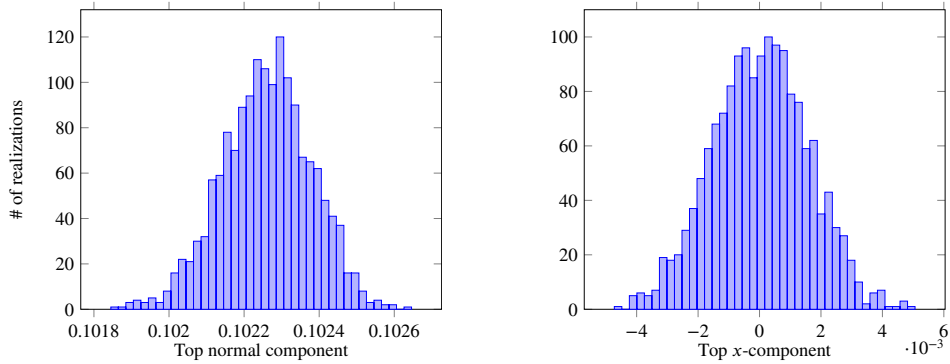

For a domain that is infinitely long in the $z$-direction, this
distribution should correspond to the estimates provided by the
homogeneized expression of
Equation~\eqref{eq:estimated-total-pressure-force}.
By constrast, the force exerted on the top and bottom faces is not
taken into account by
Equation~\eqref{eq:estimated-total-pressure-force}. For an infinte
domain, with Neumann boundary conditions on the lateral faces, this
force should be zero. However, in this case we have a finite domain,
and we impose Dirichlet boundary
conditions. Figure~\ref{fig:statistical_distribution_3d_z} shows the
statistical distribution of the normal and $x$-components on the top
face of the domain.

\subsection{A three-dimensional vascular network}\label{ssec:3dreal}

In this section we study the effect of realistic pressurized vessel
networks on a tissue sample mimicking a realistic setting of the inversion procedure
in liver elastography.
Namey, we consider a cubic tissue sample $[0,L]^3$, with $L=3$ mm
(of the order of voxel resolution of MRI scans), and elastic
characteristics similar to those found in human liver 
($\mu=2kPa$, $\lambda=50kPa$). In these settings, we investigate 
three examples of vessel distributions, with a fixed volume
fraction of $5\%$. We impose a physiological
pressure of $1kPa$.

In order to produce a realistic vessel distribution in silico, we begin with the
assumption that a vascular tree should fulfill the perfusion task with
the minimum effort, while maintaining its anatomical structure. In
general, this results in two or more competing mechanisms: on the one
hand, one expects that the total length of the vasculature shall be minimized;
on the other hand, other relevant physiological
quantity shall be minimized as well, e.g., the time needed by
oxigenated blood to reach perfusion points.

From the topological point of view, a vasculature tree can be seen as
a connected, edge-weighted undirected graph that connects some points
in the sample volume (the vertices) with a root point without any
cycles. When one tries to minimize the total edge weight, the emerging
structure is that of a minimum spanning tree\cite{Prim1957}.

In this work we use a simplified cost function, where the weight
assigned to each edge of the tree is the weighted average of two
factors: the \emph{piping cost}, represented by the Euclidean distance
between the irroration point and the connecting node in the tree, and
a \emph{total path length cost}, measuring the total path cost along
the tree from the root to the irroration point.

In particular, when the balancing factor ranges from zero to one, the
trees range from perfect minimum spanning trees (minimizing the total
length of the tree) to almost direct connections from the root to any
point (minimising the time it takes for blood to travel from the root
point to the perfusion point).

More realistic cost functions could be used\cite{Jiang2010} to take
into account other physiological details, or even mechanical
properties, but we leave this exploration for future works.  We
generate artificial vessel trees using a publicly available code\footnote{https://github.com/pherbers/MST-Dendrites} originally written to
produce synthetic neuronal structures~\cite{Cuntz2010}, setting the
balancing factor to $0.5$.

\begin{figure}[!h]
  \centering
  \includegraphics[width=.9\textwidth]{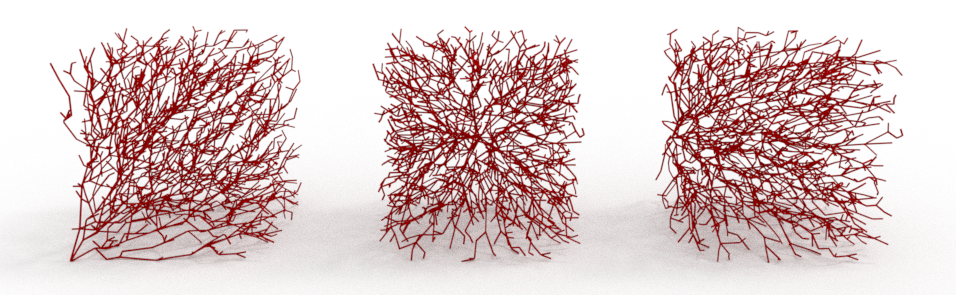}
  \caption{Randomly distributed vessels, constructed using minimum
    spanning trees with a balancing factor equal to $0.5$, irrorating
    two thousands randomly distributed points in the sample.
    Left: root point situated in the lower left corner (acronym LL).
    Center: root point located in the center of the cube (acronym C).
    Right: root point located  in the center of the left face (acronym FC).}
    \label{fig:three-vessels}
\end{figure}

We provide three different artificial vessel networks, which, in what follows,  will be denoted by LL, C, and FC, depending 
on the position of the root point relative to the sample, as
explained in Figure~\ref{fig:three-vessels}. 
For each configuration, we compute the three averaged principal
directions as the eigenvectors of the matrix
\begin{equation}
  \label{eq:principal-directions}
  \int_\Gamma \btau \otimes \btau \d\Gamma.
\end{equation}
Table~\ref{tab:stats-vessels} reports some statistical information
about the networks, providing the three principal directions and the
corresponding eigenvalues (rescaled and reordered so that the maximum
eigenvalue is always equal to one, and corresponds to the third axis
$\btau_{3}$), and the total length $L$ of the vasculature network. The
average volume fraction along the direction $\tau_i$ is proportional
to $\lambda_i$. In particular:
\begin{equation}
  \label{eq:formula-for-beta-i}
  \beta_i = \frac{L\pi a^2\lambda_i}{V_{\text{Sample}}(\lambda_1\lambda_2\lambda_3)^{\frac{1}{3}}}.
\end{equation}

\begin{table}[!h]
  \centering
  \begin{tabular}{r|r|r|r}
    & LL & C & FC \\
    \hline
    $\lambda_1, \btau_{s_1}$: &
                            $  0.4065,
                            \begin{pmatrix}
                              -0.0202 \\
                              -0.6652 \\
                              0.7463
                            \end{pmatrix}$
    &
      $        0.9264, 
      \begin{pmatrix}
        0.0775 \\-0.9225\\ -0.3780
      \end{pmatrix}
    $
    &
      $
      0.6367,
      \begin{pmatrix}
        -0.0416 \\ 0.2416  \\ 0.9694
      \end{pmatrix}
    $\\
    $\lambda_2, \btau_{s_2}$: &
                            $
                            0.4274,
                            \begin{pmatrix}
                              0.8227 \\-0.4352 \\-0.3655
                            \end{pmatrix}
    $
    &
      $
      0.9609,
      \begin{pmatrix}
        0.3186 \\  0.3822 \\ -0.8673
      \end{pmatrix}
    $
    &
      $
      0.6936,
      \begin{pmatrix}
        0.0340  \\ 0.9700 \\ -0.2403
      \end{pmatrix}
    $\\
    $\lambda_3, \btau_{s_3}$: &
                            $
                            1.0,
                            \begin{pmatrix}
                              0.5680 \\  0.6066 \\  0.5561
                            \end{pmatrix}$
    &
      $
      1.0,
      \begin{pmatrix}
        0.9446  \\-0.0531  \\ 0.3236
      \end{pmatrix}
    $
    &
      $
      1.0,
      \begin{pmatrix}
        0.9985  \\ -0.0229  \\ 0.0486
      \end{pmatrix}
    $\\[1cm]
    &    \includegraphics[height=.1\textwidth]{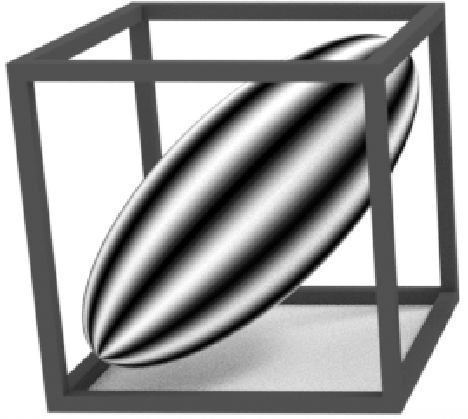}
         &  \includegraphics[height=.1\textwidth]{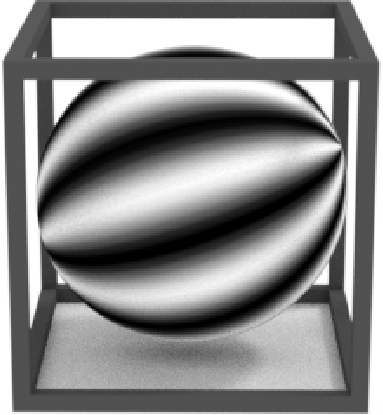}
             &  \includegraphics[height=.1\textwidth]{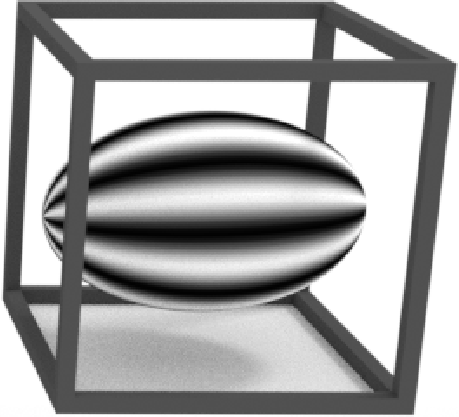}
\\
    \hline
    Total length: & $0.3839$m & $0.3795$m & $0.3878$m\\
    Vessels radius: & $3.3455\times 10^{-5}$m & $3.3646\times 10^{-5}$m & $3.3284\times 10^{-5}$m
  \end{tabular}
\caption{Statistical data of the vessels trees constructed in
  Figure~\ref{fig:three-vessels}.}
  \label{tab:stats-vessels}
\end{table}

In the LL and FC cases, alignment is predominant along one direction
(the principal diagonal in the LL case, and the $x$ axis in the FC
case), while in the C case, the alignment is roughly uniformly
distributed. Only the LL case should show a significant amount of
pressure induced shear (see
Figure~\ref{fig:pressure-induced-shear}). In all other cases, a
pressurisation would produce uniform deformations in the C case, and
non-uniform dilations along the $x$ axis and on the $yz$ plane in the
FC case, without significant induced shear.


\begin{figure}
  \centering
  \includegraphics[width=.27\textwidth,trim=0cm 0cm 14.5cm 0cm,clip=true]{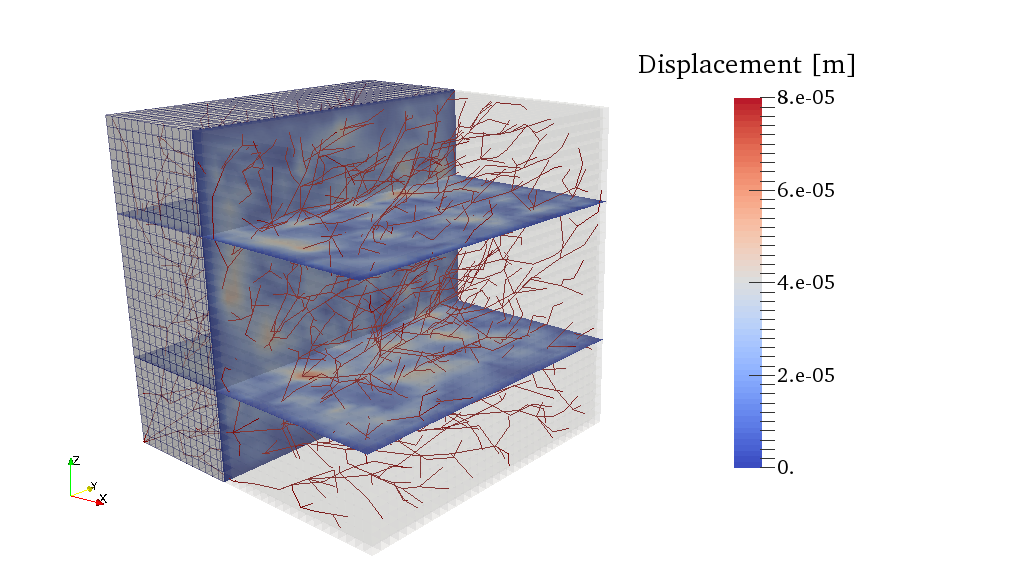}
  \includegraphics[width=.27\textwidth,trim=0cm 0cm 14.5cm 0cm,clip=true]{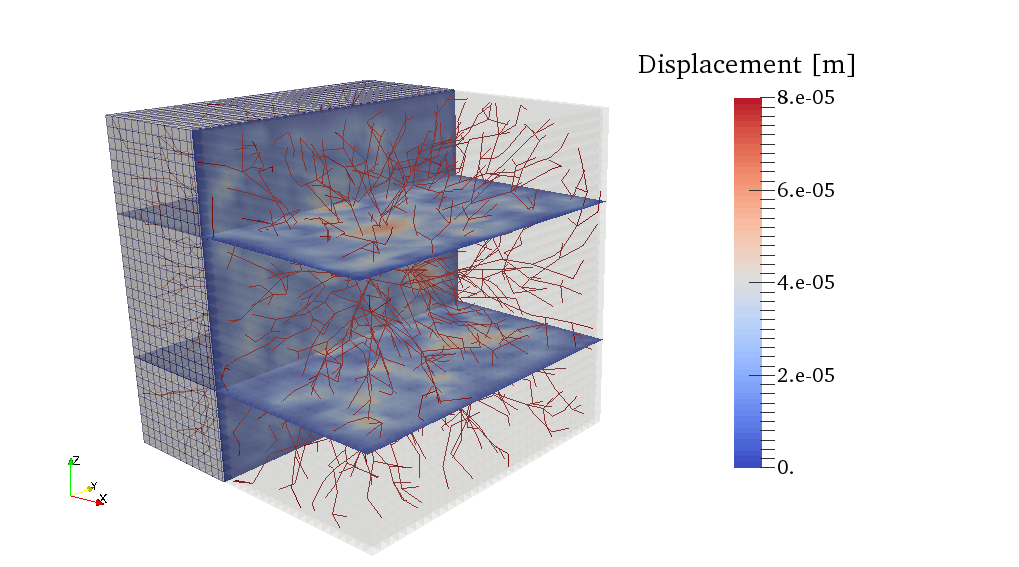}
  \includegraphics[width=.27\textwidth,trim=0cm 0cm 14.5cm 0cm,clip=true]{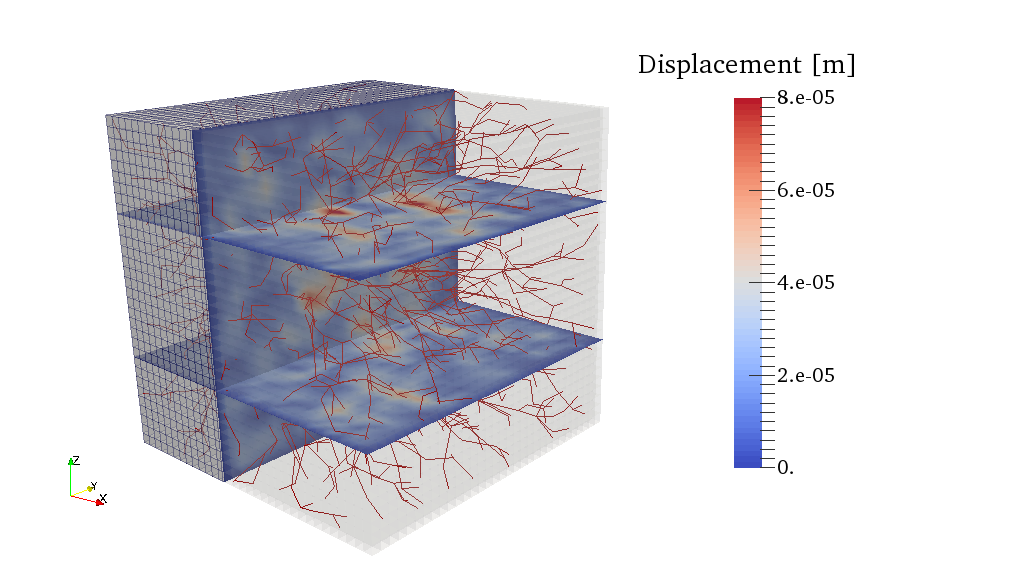}
  \includegraphics[width=.13\textwidth]{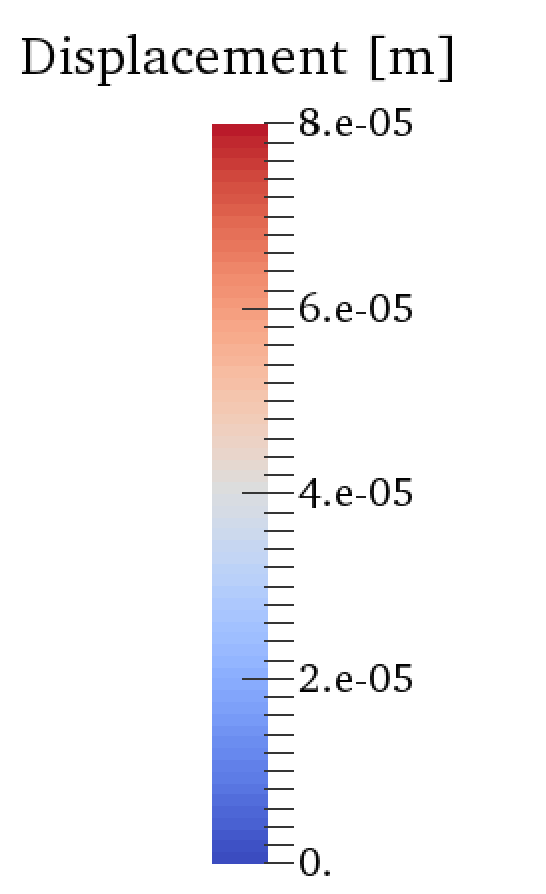}
  \caption{Magnitude of the displacement field on selected slices
in the three   sample networks (in red). The maximum displacements are given by
    $7.9\times10^{-5}$ m in the LL case (left), $9.0\times10^{-5}$ m  in the C
    case (center), and  $9.6\times10^{-5} $m  in the FC case (right).}
  \label{fig:displacement-three-vessels}
\end{figure}

For all these simulations, we constructed a structured mesh with grid
size roughly equal to $h=10^{-4}$. Given the relatively small number
of vessels in each principal direction, we do not expect the
homogeneized model presented in Section~\ref{ssec:derivation-3d} to
provide the same answers of the multiscale simulation. We measure the
pressure induced tractions defined by
\begin{equation}
  \label{eq:pressure-induced-moments}
  \mathbf{M}_{ij} = \mathbf{F}_{i+} \cdot \mathbf{e}_j -
  \mathbf{F}_{i-}  \cdot \mathbf{e}_j,
\end{equation}
where we indicate with $\mathbf{F}_{i\pm}$ the average force measured
on the face with normals $\mathbf{e}_i$ and $-\mathbf{e}_i$
respectively. These forces are the ones that an instrument would
measure on the given tissue sample. 
Table~\ref{tab:pressure-induced-moments} report these measurements.


\definecolor{verylow}{HTML}{D5F5E3} 
\definecolor{low}{HTML}{82E0AA}
\definecolor{medium}{HTML}{2ECC71} 
\definecolor{high}{HTML}{28B463} 
\definecolor{veryhigh}{HTML}{239B56} 

\begin{table}
  \centering
  \resizebox{\textwidth}{!}{
    \begin{tabular}{ccc}
       LL & C & FC \\
      $ \begin{pmatrix}
        \cellcolor{high} -3.00e-03 & \cellcolor{low} 3.30e-04 & \cellcolor{low} 2.00e-04\\
       \cellcolor{low}  3.04e-04 & \cellcolor{high} -2.84e-03 & \cellcolor{low}  4.07e-04\\
        \cellcolor{low}  2.05e-04 & \cellcolor{low} 2.12e-04 & \cellcolor{high} -2.60e-03\\
      \end{pmatrix} $ 
          & 
            $ \begin{pmatrix}
              \cellcolor{veryhigh} -5.06e-03 & \cellcolor{low} -1.33e-04 & \cellcolor{verylow} 7.93e-07\\
             \cellcolor{low} 1.28e-04 & \cellcolor{veryhigh} -5.05e-03 & \cellcolor{verylow} -6.61e-05\\
             \cellcolor{verylow} 6.16e-05 & \cellcolor{verylow} -7.52e-05 & \cellcolor{veryhigh}  -5.43e-03\\
            \end{pmatrix} $ 
          & 
            $ \begin{pmatrix}
              \cellcolor{medium}  -1.11e-03 & \cellcolor{verylow} 4.52e-06 & \cellcolor{verylow}  6.23e-05\\
              \cellcolor{verylow}  5.76e-05 &  \cellcolor{veryhigh} -4.61e-03 & \cellcolor{verylow}  1.71e-05\\
             \cellcolor{verylow}  2.25e-05 & \cellcolor{verylow} -5.25e-05 &\cellcolor{veryhigh}  -4.29e-03\\
            \end{pmatrix} $ 
    \end{tabular}
  }
  \caption{Pressure induced tractions (measured in Pa) in the three
    simulations (see Equation~\eqref{eq:pressure-induced-moments} for
    a definition). Notice how in the case of perfectly symmetric
    domains and vessels, these matrices should be symmetric, and close
    to $|A|\tens{\sigma}$. The deviation from $\tens{\sigma}$ is an
    effect of the lack of symmetry in the vessels. The color codes indicates the order of magnitude
    of the different entries of the table (increasing, from light to dark green).}
  \label{tab:pressure-induced-moments}
\end{table}

In all cases, nonzero off-diagonal terms indicate the presence of
pressure induced shear. Its magnitude and direction is non-trivially
connected with the local distribution of the vasculature. For example,
in the LL case, on each pair of opposite faces (a line in the matrices
of Table~\ref{tab:pressure-induced-moments}), we measure positive
forces (of roughly the same order of magnitude, ranging from 0.2 to 0.4 mPa) in the tangential direction. In the normal direction, instead,
we always measure a negative traction (of roughly the same intensity),
in all three cases, indicating the presence of vessels whose
directions have nonzero scalar product with the normal to \emph{every}
face.

In the C case, for every face there are always two directions along
which contraction is happening, while the FC case shows a mixture of
the LL and C behaviour (sometimes contraction along two directions,
sometimes only along the normal). From these results, it is evident
that drawing any type of conclusions on the effective mechanical 
properties based only on average measurements
on the faces may be too inaccurate, if one does not properly take into
account the effect of the vasculature.

\section{Conclusions}
\label{sec:conclusions}
We presented a multiscale modeling of biphasic tissues composed by an
elastic matrix and a set of thin vessels. 
 Our approach is based on an immersed method, treating the interface between
solid and fluid as an immersed one dimensional manifold (described by
centerline, radius, and pressure).
Neglecting long range interaction between vessels, as well as vessel curvature,
we derived a finite element formulation in which the effect of the vascular network is modeled as a singular term -- with support on the
one-dimensional vessels -- in the elasticity equations.

The immersed method allows to simulate the effective behavior of a vascularized
tissue without requiring the full resolution of the vasculature. Hence, 
the main advantage of our approach is that it drastically reduces the
computational complexity when dealing with large networks of vessels,
as it does not require to fully resolve the fluid-solid interface
within the computational mesh. In fact, vessels are represented 
by point-wise singular or hypersingular sources (regularized Dirac delta functions) distributed along the 
centerline of the vessels, and whose intensities depend on the physical and geometrical
parameters.

In order to validate the method, we show several numerical tests in simplified cases (a single vessel or
a bifurcation), comparing
the results of the immersed methods with the exact (when available) and fully resolved solutions.
In particular, we showed optimal convergence of the numerical solution in $L^2$ and $H^1$ norms
at the level of coarse scale, i.e.,  as long as the mesh resolution does not go beyond the spatial scale of the fluid vessels.
The numerical results show as well that, at the level of coarse scales, the immersed method
delivers an accuracy comparable to the fully resolved finite element solution, when using an overall comparable amount of degrees
of freedom, with the difference that the immersed formulation does not require the discretization of the 
vessel boundaries. In order to validate the method in presence of curved (with small curvature) and intersecting vessels, 
we compared the immersed solution with a full finite element solution on a mesh resolving in detail a vessel junction.

Starting from the variational formulation with immersed singular terms, we
propose a theoretical multiscale framework to analytically describe, in the case of uniform vessel distribution, 
the effect of a pressurized vasculature onto the mechanical properties of the resulting tissue. 
The purpose of these analysis is to, firstly, set the basis for the understanding on how vascular structures
should be taken into account when shear and compression experiments are used to characterize tissues, and, secondly,
to provide an explanation for experimental observations revealing strong differences in the
mechanical parameters (estimated, e.g.,  using elastography) in vivo and ex vivo.

The outcome of the multiscale analysis has been validated with a set
of statistical simulations, in order to correlate effective material
properties with the volume fraction of the vasculature and vessel
pressure.  The presented examples showed that the non-matching
immersed method can be used to investigate, in silico and from a
statistical perspective, the mechanical behavior of the tissue given
the (statistical) properties of the vasculature.  To this respect, the
main advantage of the immersed formulation is that the same
computational mesh can be employed for simulating different
realization, since the vessels are not explicitly resolved by the
mesh.

In order to demonstrate the potential of the immersed method in
realistic cases, we simulate a cubic tissue sample with three
different realization of a randomly generated vasculature, prescribing
different geometrical properties, and yielding different degrees of
anisotropy.
The scope of these simulations is to estimate, in silico, shear and
compression effects due to the presence of pressurised vessel
networks. This information provides a tools to correct the estimation
of mechanical parameters of in-vivo tissues.
We show that the pressure in the 
(microscopic) vasculatures induces a (macroscopic) shear, measured as the difference between forces acting 
on boundary faces.  Although a validation
with a fully resolved simulation is not possible in this case, our results shows that the pressure induced shear
is highly correlated with the orientation of the underlying vasculature.

These results shall be interpreted as an intermediate step towards the
estimation of mechanical parameters of vascularized tissues. Our
purpose is to employ the proposed formulation for the solution of
inverse problems targeting the characterization of effective tissue
properties.  The homogenized characterization, the statistical
simulations, and the simulation considering randomly generated
vasculature aims at providing evidence that the immersed method can be
effectively used to link microscale properties (geometry and pressure
of blood vessels) with macroscopic parameters (shear and compression
modulus of a tissue sample). In this work, we presented a detailed
numerical validation, limited to qualitative aspects in the case of
complex vasculatures. In order to assess the full potential of the
immersed method to address inverse problems in relevant clinical
contexts, future work shall focus -- in collaboration with
experimentalists -- on a more detailed experimental validation.


This work is motivated by the non-invasive estimation of mechanical
properties of living tissues using magnetic resonance elastography
(MRE), which combines displacement fields acquired using phase
contrast MRI (whose typical voxel resolution, in this context, is of
the order of millimeters) with a mathematical model of tissue
mechanics.
Physical models currently used in elastography are mainly restricted
to isotropic and homogeneous elastic or viscoelastic tissues. In
particular, due to the need of keeping the computational cost low,
clinical applications are often based on linear elasticity.
One of the scopes of our work is to go beyond these models, taking into
account sources of non linearity and anisotropy due to an underlying
vasculature, without necessarily sacrificing the computational
efficiency of the method.
The multiscale model derived in this first paper is based on the assumption of that the tissue matrix 
behaves as a linear elastic isotropic tissue. However, by taking into account the fluid vasculature, 
non linearities and anisotropies arise at the effective tissue scale.

On the other hand, although the assumption of linear behavior can be
justified in the range of small displacement induced in MRE, we plan
to extend the proposed model to more general structural models (e.g.,
poroelasticity) in upcoming works.

The presented multiscale model is limited to the situation of a static
fluid (at a given pressure) and it does not consider a full
fluid-structure interaction problem, i.e., we neglect the effect of
the tissue pressure onto the vascular pressure, and we ignore long
range interaction between far sections of the vessels.  These
assumptions are motivated by the fact that, in the context of MRE, the
frequency of induced harmonic excitation (30 to 60Hz) is much higher
than a typical heart beat, and fluid pressure can therefore considered
to be constant.
One of the goal of our model is the possibility of linking
effective mechanical parameters with possibly pathological hemodynamic, motivated by current research
in obtaining non-invasive biomarkers of hypertension via MRE\cite{hirsch-etal-2013-compmre,hirsch-etal-2014-liver}.
Therefore, taking into account the coupling with an active vasculature, not only limited to the vessels irrigating the tissue sample,
is an aspect of utmost importance, and the coupling of the proposed immersed method with a
time dependent one-dimensional blood flow model\cite{mueller-etal-16} is subject of ongoing research.

In this context, it shall be also observed that the proposed finite element method
can be used to derive efficient reduced order models (e.g., reduced basis method or 
proper orthogonal decomposition), where the finite element matrix
is assembled only once and projected onto a small subspace.
in fact, since the vasculature enter the elasticity equations only as a singular right hand side,
variation in the vessel (e.g., pressure depending on time) do not require reassembling the linear system, and all operations
can be performed within the reduced space.

\bibliographystyle{abbrv}
\bibliography{ref}

\end{document}

%% file: tables/error_exact_hole_a.1.tex
\begin {tabular}{|r<{\pgfplotstableresetcolortbloverhangright }@{}l<{\pgfplotstableresetcolortbloverhangleft }|r<{\pgfplotstableresetcolortbloverhangright }@{}l<{\pgfplotstableresetcolortbloverhangleft }|r<{\pgfplotstableresetcolortbloverhangright }@{}l<{\pgfplotstableresetcolortbloverhangleft }c|r<{\pgfplotstableresetcolortbloverhangright }@{}l<{\pgfplotstableresetcolortbloverhangleft }c|}%
\hline \multicolumn {2}{|r}{\#cells}&\multicolumn {2}{|r|}{\#dofs}&\multicolumn {2}{c}{$\|u-u_h\|_{L^2}$}&&\multicolumn {2}{c}{$\|u-u_h\|_{H^1}$}&\\\hline \hline %
\rowcolor [gray]{.95}$16$&$$&$48$&$$&$5$&$.5160\cdot 10^{-3}$&$-$&$1$&$.1130\cdot 10^{-1}$&$-$\\%
$64$&$$&$160$&$$&$2$&$.8540\cdot 10^{-3}$&0.95&$7$&$.7860\cdot 10^{-2}$&0.51\\%
\rowcolor [gray]{.95}$256$&$$&$576$&$$&$1$&$.2270\cdot 10^{-3}$&1.22&$4$&$.8900\cdot 10^{-2}$&0.67\\%
$1{,}024$&$$&$2{,}176$&$$&$4$&$.0740\cdot 10^{-4}$&1.59&$2$&$.7400\cdot 10^{-2}$&0.84\\%
\rowcolor [gray]{.95}$4{,}096$&$$&$8{,}448$&$$&$1$&$.1330\cdot 10^{-4}$&1.85&$1$&$.4290\cdot 10^{-2}$&0.94\\%
$16{,}384$&$$&$33{,}280$&$$&$2$&$.9240\cdot 10^{-5}$&1.95&$7$&$.2320\cdot 10^{-3}$&0.98\\%
\rowcolor [gray]{.95}$65{,}536$&$$&$132{,}096$&$$&$7$&$.3720\cdot 10^{-6}$&1.99&$3$&$.6280\cdot 10^{-3}$&1.00\\%
$262{,}144$&$$&$526{,}336$&$$&$1$&$.8470\cdot 10^{-6}$&2.00&$1$&$.8150\cdot 10^{-3}$&1.00\\\hline %
\end {tabular}%

%% file: tables/error_approx_dirac_a.1.tex
\begin {tabular}{|r<{\pgfplotstableresetcolortbloverhangright }@{}l<{\pgfplotstableresetcolortbloverhangleft }|r<{\pgfplotstableresetcolortbloverhangright }@{}l<{\pgfplotstableresetcolortbloverhangleft }|r<{\pgfplotstableresetcolortbloverhangright }@{}l<{\pgfplotstableresetcolortbloverhangleft }c|r<{\pgfplotstableresetcolortbloverhangright }@{}l<{\pgfplotstableresetcolortbloverhangleft }c|}%
\hline \multicolumn {2}{|r}{\#cells}&\multicolumn {2}{|r|}{\#dofs}&\multicolumn {2}{c}{$\|u-u_h\|_{L^2}$}&&\multicolumn {2}{c}{$\|u-u_h\|_{H^1}$}&\\\hline \hline %
\rowcolor [gray]{.95}$80$&$$&$178$&$$&$1$&$.0750\cdot 10^{-2}$&$-$&$1$&$.6380\cdot 10^{-1}$&$-$\\%
$320$&$$&$674$&$$&$4$&$.1020\cdot 10^{-3}$&1.39&$1$&$.2130\cdot 10^{-1}$&0.43\\%
\rowcolor [gray]{.95}$1{,}280$&$$&$2{,}626$&$$&$7$&$.3670\cdot 10^{-4}$&2.48&$3$&$.0780\cdot 10^{-2}$&1.98\\%
$5{,}120$&$$&$10{,}370$&$$&$4$&$.8270\cdot 10^{-4}$&0.61&$1$&$.5250\cdot 10^{-2}$&1.01\\%
\rowcolor [gray]{.95}$20{,}480$&$$&$41{,}218$&$$&$4$&$.4570\cdot 10^{-4}$&0.12&$7$&$.6150\cdot 10^{-3}$&1.00\\%
$81{,}920$&$$&$164{,}354$&$$&$4$&$.4280\cdot 10^{-4}$&0.01&$4$&$.8220\cdot 10^{-3}$&0.66\\%
\rowcolor [gray]{.95}$327{,}680$&$$&$656{,}386$&$$&$4$&$.4250\cdot 10^{-4}$&0.00&$3$&$.9090\cdot 10^{-3}$&0.30\\\hline %
\end {tabular}%

%% file: tables/error_exact_hole_a.01.tex
\begin {tabular}{|r<{\pgfplotstableresetcolortbloverhangright }@{}l<{\pgfplotstableresetcolortbloverhangleft }|r<{\pgfplotstableresetcolortbloverhangright }@{}l<{\pgfplotstableresetcolortbloverhangleft }|r<{\pgfplotstableresetcolortbloverhangright }@{}l<{\pgfplotstableresetcolortbloverhangleft }c|r<{\pgfplotstableresetcolortbloverhangright }@{}l<{\pgfplotstableresetcolortbloverhangleft }c|}%
\hline \multicolumn {2}{|r}{\#cells}&\multicolumn {2}{|r|}{\#dofs}&\multicolumn {2}{c}{$\|u-u_h\|_{L^2}$}&&\multicolumn {2}{c}{$\|u-u_h\|_{H^1}$}&\\\hline \hline %
\rowcolor [gray]{.95}$16$&$$&$48$&$$&$5$&$.0980\cdot 10^{-4}$&$-$&$7$&$.0660\cdot 10^{-3}$&$-$\\%
$64$&$$&$160$&$$&$2$&$.6450\cdot 10^{-4}$&0.95&$8$&$.5960\cdot 10^{-3}$&-0.28\\%
\rowcolor [gray]{.95}$256$&$$&$576$&$$&$1$&$.3980\cdot 10^{-4}$&0.92&$1$&$.0270\cdot 10^{-2}$&-0.26\\%
$1{,}024$&$$&$2{,}176$&$$&$9$&$.2200\cdot 10^{-5}$&0.60&$1$&$.0350\cdot 10^{-2}$&-0.01\\%
\rowcolor [gray]{.95}$4{,}096$&$$&$8{,}448$&$$&$6$&$.1250\cdot 10^{-5}$&0.59&$8$&$.5820\cdot 10^{-3}$&0.27\\%
$16{,}384$&$$&$33{,}280$&$$&$3$&$.1920\cdot 10^{-5}$&0.94&$5$&$.9620\cdot 10^{-3}$&0.53\\%
\rowcolor [gray]{.95}$65{,}536$&$$&$132{,}096$&$$&$1$&$.2130\cdot 10^{-5}$&1.40&$3$&$.5610\cdot 10^{-3}$&0.74\\%
$262{,}144$&$$&$526{,}336$&$$&$3$&$.6130\cdot 10^{-6}$&1.75&$1$&$.9160\cdot 10^{-3}$&0.89\\\hline %
\end {tabular}%

%% file: tables/error_approx_dirac_a.01.tex
\begin {tabular}{|r<{\pgfplotstableresetcolortbloverhangright }@{}l<{\pgfplotstableresetcolortbloverhangleft }|r<{\pgfplotstableresetcolortbloverhangright }@{}l<{\pgfplotstableresetcolortbloverhangleft }|r<{\pgfplotstableresetcolortbloverhangright }@{}l<{\pgfplotstableresetcolortbloverhangleft }c|r<{\pgfplotstableresetcolortbloverhangright }@{}l<{\pgfplotstableresetcolortbloverhangleft }c|}%
\hline \multicolumn {2}{|r}{\#cells}&\multicolumn {2}{|r|}{\#dofs}&\multicolumn {2}{c}{$\|u-u_h\|_{L^2}$}&&\multicolumn {2}{c}{$\|u-u_h\|_{H^1}$}&\\\hline \hline %
\rowcolor [gray]{.95}$80$&$$&$178$&$$&$3$&$.2390\cdot 10^{-4}$&$-$&$2$&$.5120\cdot 10^{-2}$&$-$\\%
$320$&$$&$674$&$$&$2$&$.0680\cdot 10^{-4}$&0.65&$9$&$.8780\cdot 10^{-3}$&1.35\\%
\rowcolor [gray]{.95}$1{,}280$&$$&$2{,}626$&$$&$2$&$.0690\cdot 10^{-4}$&-0.00&$1$&$.9750\cdot 10^{-2}$&-1.00\\%
$5{,}120$&$$&$10{,}370$&$$&$1$&$.2670\cdot 10^{-4}$&0.71&$1$&$.5670\cdot 10^{-2}$&0.33\\%
\rowcolor [gray]{.95}$20{,}480$&$$&$41{,}218$&$$&$6$&$.4070\cdot 10^{-5}$&0.98&$1$&$.3940\cdot 10^{-2}$&0.17\\%
$81{,}920$&$$&$164{,}354$&$$&$1$&$.6500\cdot 10^{-5}$&1.96&$4$&$.7070\cdot 10^{-3}$&1.57\\%
\rowcolor [gray]{.95}$327{,}680$&$$&$656{,}386$&$$&$1$&$.3990\cdot 10^{-5}$&0.24&$2$&$.0450\cdot 10^{-3}$&1.20\\\hline %
\end {tabular}%